\documentclass[12pt]{article}
\usepackage[american]{babel}

\let\LatexGamma\Gamma
\usepackage{amsmath}
\usepackage{latexsym}
\usepackage{amssymb}
\usepackage{theorem}
\usepackage[arrow,matrix,curve,poly]{xy}
\theoremstyle{change}
\newtheorem{Thm}{Theorem}[subsection]
\newtheorem{Cor}[Thm]{Corollary}
\newtheorem{Prop}[Thm]{Proposition}
\newtheorem{Lem}[Thm]{Lemma}
{\theorembodyfont{\rmfamily}
\newtheorem{Num}[Thm]{}

}
\renewcommand{\Form}{\mathsf{Form}}
\newcommand{\Alt}{\mathsf{Alt}}
\newcommand{\Quad}{\mathsf{Quad}}
\newcommand{\Herm}{\mathsf{Herm}}
\newcommand{\End}{\mathsf{End}}
\newcommand{\Hom}{\mathsf{Hom}}

\newcommand{\PG}{\mathsf{PG}}
\newcommand{\Cen}{\mathsf{Cen}}
\newcommand{\Int}{\mathsf{Int}}
\newcommand{\Aut}{\mathsf{Aut}}
\newcommand{\Out}{\mathsf{Out}}
\newcommand{\proof}{\par\medskip\rm\emph{Proof. }}
\newcommand{\qed}{\ \hglue 0pt plus 1filll $\Box$}
\renewcommand{\SS}{\mathbb{S}}
\newcommand{\RR}{\mathbb{R}}
\newcommand{\ZZ}{\mathbb{Z}}
\newcommand{\FF}{\mathbb{F}}
\newcommand{\HH}{\mathbb{H}}
\newcommand{\sK}{\mathsf{K}}
\newcommand{\cL}{\mathcal{L}}
\newcommand{\cP}{\mathcal{P}}

\newcommand{\id}{\mathsf{id}}
\newcommand{\Sp}{\mathsf{Sp}}
\newcommand{\SL}{\mathsf{SL}}
\newcommand{\U}{\mathsf{U}}
\newcommand{\GL}{\mathsf{GL}}
\newcommand{\EL}{\mathsf{EL}}
\newcommand{\PGL}{\mathsf{PGL}}
\newcommand{\PEL}{\mathsf{PEL}}
\newcommand{\PSL}{\mathsf{PSL}}
\newcommand{\GammaL}{\mathsf{\LatexGamma L}}
\newcommand{\PGammaL}{\mathsf{P\LatexGamma L}}
\newcommand{\GammaU}{\mathsf{\LatexGamma U}}
\newcommand{\GU}{\mathsf{GU}}

\newcommand{\eps}{\varepsilon}
\renewcommand{\emptyset}{\varnothing}
\newcommand{\Gr}{\mathsf{Gr}}
\renewcommand{\dim}{{\mathsf{dim}}}
\renewcommand{\det}{{\mathsf{det}}}
\newcommand{\fin}{{\mathsf{fin}}}
\newcommand{\stable}{{\mathsf{stb}}}
\newcommand{\Sym}{\mathsf{Sym}}
\newcommand{\op}{\mathsf{op}}
\newcommand{\rk}{\mathsf{rk}}
\renewcommand{\max}{\mathsf{max}}
\newcommand{\ind}{\mathsf{ind}}
\newcommand{\Mod}{\mathfrak{Mod}}
\newcommand{\kernel}{\mathsf{ker}}
\newcommand{\bra}[1]{\langle#1\rangle}
\begin{document}

\title{\bf Buildings and Classical Groups}
\author{Linus Kramer\thanks{Supported by a Heisenberg fellowship
by the Deutsche Forschungsgemeinschaft}\\
\small Mathematisches Institut,
Universit\"at W\"urzburg\\
\small
Am Hubland,
D--97074 W\"urzburg,
Germany \\
\small email: {\tt kramer@mathematik.uni-wuerzburg.de}}
\date{}
\maketitle

In these notes we describe the classical groups, that is, the linear
groups and the orthogonal, symplectic, and unitary groups, acting on
finite dimensional vector spaces over skew fields, as well as their
pseudo-quadratic generalizations. Each such group corresponds 
in a natural way to
a point-line geometry, and to a spherical building. The geometries in
question are projective spaces and polar spaces. We emphasize in
particular the r\^ole played by root elations and the groups generated
by these elations. The root elations reflect --- via their commutator
relations --- algebraic properties of the underlying vector space.

We also discuss some related algebraic topics: the classical groups as
permutation groups and the associated simple groups. I have included
some remarks on K-theory, which might be interesting for applications.
The first K-group measures the difference between the classical group
and its subgroup generated by the root elations. The second K-group
is a kind of fundamental group of the group generated by the root
elations and is related to central extensions.
I also included some material on Moufang sets, since this is an
interesting topic. In this context, the projective line over a skew
field is treated in some detail, and possibly with some new results.
The theory of unitary groups is developed along the lines of
Hahn \& O'Meara \cite{HOM}. Other important sources are the books
by Taylor \cite{Taylor} and Tits \cite{TitsLNM}, and
the classical books by Artin \cite{Artin}
and Dieudonn\'e \cite{Dieu}. The books by Knus \cite{Knus}
and W. Scharlau \cite{WScharlau} should also be mentioned here.
Finally, I would like to recommend the surveys by Cohen \cite{Cohen}
and R. Scharlau \cite{RScharlau}.

While most of these matters are well-known to experts, there seems to be
no book or survey article
which contains these aspects simultaneously. Taylor's book \cite{Taylor}
is a nice and readable introduction to classical groups, but it is clear
that the author secretly thinks of finite fields --- the non-commutative
theory is almost non-existent in his book. On the other extreme, the book
by Hahn \& O'Meara \cite{HOM}
contains many deep algebraic facts about classical
groups over skew fields; however, their book contains virtually no
geometry (the key words \emph{building} or \emph{parabolic subgroup} are
not even in the index!). So, I hope
that this survey --- which is based on lectures notes from a course I
gave in 1998 in W\"urzburg
 --- is a useful compilation of material from different
sources.

There is (at least)
one serious omission (more omissions can be found in the last
section): I have included nothing about
coordinatization. The coordinatization of these geometries is another
way to recover some algebraic structure from geometry. For projective
planes and projective
spaces, this is a classical topic, and I just mention the
books by Pickert \cite{Pi} and Hughes \& Piper \cite{HP}.
For generalized quadrangles (polar spaces of rank $2$), Van Maldeghem's
book \cite{HVM} gives a comprehensive introduction. Understanding
coordinates in the rank 2 case is in most cases sufficient in order
to draw some algebraic conclusions.

I would like to thank Theo Grundh\"ofer and Bernhard M\"uhlherr
for various remarks on the manuscript, and to Katrin Tent.
Without her, these notes would have become hyperbolic.

\tableofcontents
\SilentMatrices

\subsection*{Preliminaries}
We first fix some algebraic terminology. A (left) action of a group $G$
on a set $X$ is a homomorphism $\xymatrix@1{G\ar[r]&{\Sym(X)}}$ of $G$ 
into the permutation group $\Sym(X)$
of $X$; the permutation induced by $g$ is (in most cases) also denoted by
$g=[\xymatrix@1{x\ar@{|->}[r]&g(x)}]$.
A \emph{right action} is an anti-homomorphism
$\xymatrix@1{G\ar[r]&{\Sym(X)}}$; for right
actions, we use exponential notation,
$\xymatrix@1{x\ar@{|->}[r]&x^g}$. The set of
all $k$-element subsets of a set $X$ is denoted $\binom Xk$.

Given a (not necessarily commutative)
ring $R$ (with unit $1$), we let $R^\times$ denote the group of
multiplicatively invertible elements.
If $R^\times=R\setminus\{0\}$, then $R$ is called
a \emph{skew field} or \emph{division ring}.
Occasionally, we use Wedderburn's Theorem.

\medskip\noindent
\textbf{Wedderburn's Theorem}
\emph{A finite skew field is commutative and thus isomorphic to some
Galois field $\FF_q$, where $q$ is a prime power.}

\medskip\noindent
For a proof see Artin \cite{Artin} Ch.~I Thm.~1.14 or
Grundh\"ofer \cite{Grundh}.
\qed

\medskip\noindent
The \emph{opposite ring}
$R^\op$ of a ring $R$
is obtained by defining a new multiplication $a\cdot b=ba$
on $R$.
An \emph{anti-automorphism} $\alpha$ of $R$ is a ring isomorphism
$\xymatrix@1{R\ar[r]^\alpha&{R^\op}}$, i.e. $(xy)^\alpha=y^\alpha x^\alpha$.
The group consisting of all automorphisms and anti-automorphisms
of $D$ is denoted $\mathsf{AAut}(D)$; it has the automorphism group
$\Aut(D)$ of $D$ as a normal subgroup (of index $1$ or $2$).
Let $M$ be an abelian
group. A \emph{right $R$-module structure on $M$} is a ring homomorphism
\[
\xymatrix{{R^\op}\ar[r]^(.4)\rho&{\End(M);}}
\]
as customary,  we write 
\[
mrs=\rho(s\cdot r)(m)
\]
for
$r,s\in R$ and $m\in M$ ('scalars to the right'). The abelian
category of all right $R$-modules is defined
in the obvious way and denoted $\Mod_R$; the subcategory consisting
of all finitely generated right $R$-modules is denoted $\Mod_R^\fin$.
If $D$ is a skew field, then $\Mod_D^\fin$ is the category of all
finite dimensional \emph{vector spaces}
over $D$. Similarly, we define the category
of left $R$-modules $_R\Mod$;
given a right $R$-module $M$, we have the \emph{dual}
$M^\vee=\Hom_R(M,R)$ which is
in a natural way a left $R$-module.

A right module over a skew field will be called a \emph{right vector space}
or just a vector space. Mostly, we will consider (finite dimensional)
right vector
spaces (so linear maps act from the left and scalars act from
the right), but occasionally we will need both types. Of course, all these
distinctions are obsolete over commutative skew fields, but in the
non-commutative case one has to be careful.

\section{Projective geometry and the general linear group}

In this first part we consider the projective geometry over
a skew field $D$ and the related groups.

\subsection{The general linear group}
\begin{quote}
\emph{We introduce the projective geometry $\PG(V)$ associated to a
finite dimensional vector space $V$ and the general linear
group $\GL(V)$, as well as the general semilinear group $\GammaL(V)$.
We describe the relations between these groups and their projective
versions $\PGL(V)$ and $\PGammaL(V)$. Finally, we
recall the 'first' Fundamental Theorem of Projective Geometry.}
\end{quote}
Let $V$ be a right vector space
over a skew field $D$, of
(finite) dimension $\dim(V)=n+1\geq 2$.
The collection of all $k$-dimensional subspaces of $V$
is the \emph{Grassmannian}
\[
\Gr_k(V)=\{X\subseteq V|\ \dim(X)=k\}.
\]
The elements of $\Gr_1(V)$, $\Gr_2(V)$ and $\Gr_n(V)$ are called
\emph{points}, \emph{lines}, and \emph{hyperplanes}, respectively.
Two subspaces $X,Y$ are called \emph{incident}, 
\[
X * Y,
\]
if either
$X\subseteq Y$ or $Y\subseteq X$. The resulting $n$-sorted structure
\[
\PG(V)=(\Gr_1(V),\ldots,\Gr_n(V),*)
\]
is the \emph{projective geometry} of rank $n$ over $D$.

It is clear that every linear bijection of $V$ induces
an automorphism of $\PG(V)$. More generally, every semilinear bijection
induces an automorphism of $\PG(V)$. Recall that a group
endomorphism $f$ of $(V,+)$ is called \emph{semilinear} (relative to an
automorphism
$\theta$ of $D$) if 
\[
f(va)=f(v)a^\theta
\]
holds for all $a\in D$ and
$v\in V$. The group of all semilinear bijections of $V$ is denoted
$\GammaL(V)$; it splits as a semidirect product with the
\emph{general linear group} $\GL(V)$, the
groups consisting of all linear bijections,
as a normal subgroup,
\[
\xymatrix{
1 \ar[r]  & {\GL(V)} \ar[r] & {\GammaL(V)} \ar[r]^{\longleftarrow}
 & {\Aut(D)} \ar[r] & 1.}
\]
As usual, we write $\GL(D^n)=\GL_n(D)$, and similarly for the groups
$\PGL(V)$ and $\PGammaL(V)$ induced on $\PG(V)$.

The kernel of the action of $\GL(V)$ on $\PG(V)$ consists of all maps of the
form
$\xymatrix@1{\rho_c:v\ar@{|->}[r]&vc}$, where $c\in\Cen(D^\times)$.
Similarly, the kernel of
the action of $\GammaL(V)$ consists of all maps of the form
$\xymatrix@1{\rho_c:v\ar@{|->}[r]&vc}$,
for $c\in D^\times$. These groups fit together in a
commutative diagram with exact rows and columns
\[
\xymatrix{
& 1 \ar[d] & 1 \ar[d] & 1 \ar[d] \\
1 \ar[r] & {\Cen(D^\times)} \ar[r]\ar[d] & {D^\times}       \ar[r]\ar[d]
& {\Int(D)} \ar[r]\ar[d] & 1 \\
1 \ar[r] & {\GL(V)} \ar[r]\ar[d] & {\GammaL(V)} \ar[r]\ar[d] 
& {\Aut(D)} \ar[r]\ar[d] & 1 \\
1 \ar[r] & {\PGL(V)} \ar[r]\ar[d] & {\PGammaL(V)} \ar[r]\ar[d] 
& {\Out(D)} \ar[r]\ar[d] & 1  \\
& 1 & 1 & 1 }
\]
as is easily checked, see Artin \cite{Artin} II.10, p.~93.
Here $\Int(D)$ is the group of inner automorphisms
$\xymatrix@1{d\ar@{|->}[r]&d^a=a^{-1}da}$ of $D$,
where $a\in D^\times$, and $\Out(D)$
is the quotient group $\Aut(D)/\Int(D)$.
The bottom line shows the groups induced on $\PG(V)$, and
the top line the kernels of the respective actions.
For example, $\Out(\HH)=1$ holds for the quaternion division algebra
$\HH$ over any real closed field, whence
$\PGL_{n+1}(\HH)\cong\PGammaL_{n+1}(\HH)$
for $n\geq 1$. On the other hand, $\Int(D)=1$ if $D$ is commutative.
Recall the Fundamental Theorem of Projective Geometry,
which basically says that $\Aut(\PG(V))=\PGammaL(V)$.
The actual statement is in fact somewhat stronger.

\begin{Num}\textbf{The Fundamental Theorem of Projective Geometry, I}
\label{FTPG1}\ 

\noindent\em
For $i=1,2$, let $V_i$ be vector spaces over $D_i$, of (finite) dimensions
$n_i+1\geq 3$, and let
\[
\xymatrix{{\PG(V_1)}\ar[r]^\phi_\cong &{\PG(V_2)}}
\]
be an isomorphism. Then $n_1=n_2$, and
there exists a skew field isomorphism
$\xymatrix@1{{\theta:D_1}\ar[r]^(.6)\cong &{D_2}}$
and a $\theta$-semilinear bijection 
\[
\xymatrix{{f:V_1}\ar[r] & V_2}
\]
(i.e. $f(va)=f(v)a^\theta$)) such that $\phi(X)=f(X)$ for all
$X\in\Gr_k(V_1)$, for $k=1,2,\ldots,n_1$. In particular,
\[
\Aut(\PG(V))=\PGammaL(V).
\]
\emph{For a proof see Hahn \& O'Meara \cite{HOM} 3.1.C,
Artin \cite{Artin} Ch.~II Thm.~2.26, L\"uneburg
\cite{Lun1} \cite{Lun2}, or Faure \& Fr\"olicher \cite{FauFro1}.}
\qed
\end{Num}
\emph{The theorem is also valid for infinite dimensional projective spaces;
this will be important in the second part, when we consider hermitian
forms.}
Note that for $n=1$, the structure $\PG(V)=(\Gr_1(V),=)$
is rather trivial. We will come back to a refined
version of the projective line in Section \ref{ProjectiveLine}.

\subsection{Elations, transvections, and the elementary linear group}
\begin{quote}
\emph{We introduce elations, transvections,
and the elementary linear group $\EL(V)$
generated by all elations. The projective elementary group is
a simple group, except for some low dimensional cases over
fields of small cardinality;
in the commutative case, the elementary linear group coincides
with the special linear group $\SL(V)$.}
\end{quote}
As in the previous section, $V$ is a right vector space over $D$ of
finite dimension $n+1\geq 2$.
Let $A\in\Gr_n(V)$ be a hyperplane, and let $z\in\Gr_1(V)$ be a point
incident with $A$. An automorphism $\tau$ of $\PG(V)$ which fixes $A$ pointwise
and $z$ linewise is called an \emph{elation} or \emph{translation},
with \emph{axis} $A$ and
\emph{center} $z$.
We can choose a base $b_0,\ldots,b_n$ of $V$ such that $A$ is spanned
by $b_0,\ldots,b_{n-1}$ and $z$ by $b_0$, and such that $\tau$ is
represented by a matrix of the form
\[
\begin{pmatrix}
1&&&&& a\\
 & 1 \\
 && 1 \\
 &&& 1 \\
 &&&& \ddots \\
 &&&&& 1
\end{pmatrix} 
\]
for some $a\in D$. Such an elation is called a $(z,A)$-elation; the
group $U_{(z,A)}$ consisting of all $(z,A)$-elations is isomorphic to 
the additive group $(D,+)$.
There is a coordinate-free way to describe elations.
Let 
\[
V^\vee=\Hom_D(V,D)
\]
denote the dual of $V$. This is a left
vector space over $D$, so we
can form the tensor product $V\otimes_D V^\vee$; this is an abelian
group (even a ring) which is naturally isomorphic to the endomorphism
ring of $V$,
\[
\End_D(V)\cong V\otimes_D V^\vee.
\]
We write
$x\otimes\phi=x\phi=[\xymatrix@1{v\ar@{|->}[r]&x\phi(v)}]$. If the elements
of $V$ are represented as column vectors and the elements of
$V^\vee$ as row vectors,
then $x\phi$ is just the standard matrix product.
Now let $(u,\rho)\in V\times V^\vee$ be a pair such that the
endomorphism $u\rho\in\End_D(V)$ is nilpotent, i.e.
$u\rho u\rho=u\rho(u)\rho=0$
(which is equivalent to $\rho(u)=0$). The linear map
\[
\tau_{u\rho}=\id_V+u\rho=[\xymatrix{v\ar@{|->}[r]&v+u\rho(v)}]
\]
is called a \emph{transvection.} Note that
\[
\tau_{ua\rho}\tau_{ub\rho}=\tau_{u(a+b)\rho}
\qquad\text{ whence }\qquad
\tau_{u\rho}\tau_{-u\rho}=\id_V.
\]
Such a transvection induces an elation
on $\PG(V)$, and conversely, every elation in $\PG(V)$ is induced by
a transvection. In fact, suppose that $u\rho\neq 0$. The center
of $\tau_{u\rho}$ is $z=uD\in\Gr_1(V)$, and the axis is
$A=\kernel(\rho)$. The group $(D,+)\cong\{\tau_{ua\rho}|\ a\in D\}$ maps
isomorphically onto $U_{(z,A)}$.
Suppose that $\tau_{u\phi}\neq 1\neq \tau_{v\psi}$ are commuting transvections,
\[
\tau_{u\phi}\tau_{v\psi}=\tau_{v\psi}\tau_{u\phi}.
\]
This implies that $u\phi(v)\psi=v\psi(u)\phi$, and thus
$\psi(u)=\phi(v)=0$.
\begin{Lem}
Two non-trivial transvections commute if and only if they have either the
same axis or the same center.
\qed
\end{Lem}

The group generated by all transvections is the
\emph{elementary linear group} $\EL(V)$. It is a normal subgroup of
the group $\GammaL(V)$ and of $\GL(V)$
(because the conjugate of a transvection is again a transvection).
The subgroup in $\Aut(\PG(V))$ generated by all elations is the
\emph{little projective group} of $\PG(V)$; it is normal and
an epimorphic image of $\EL(V)$. We denote this group by $\PEL(V)$.
We collect a few facts about the group $\EL(V)$.
\begin{Lem}
The action of $\PEL(V)$ on $\Gr_1(V)$ is 2-transitive and in particular
primitive.

\proof
It suffices to show that given a point $p\in\Gr_1(V)$, the stabilizer
$\PEL(V)_p$ acts transitively on $\Gr_1(V)\setminus\{p\}$. Also, it is
easy to see that given a line $\ell$ passing through $p$, the stabilizer
$\PEL(V)_{p,\ell}$ acts transitively on the points lying on $\ell$ and
different from $p$ --- choose an axis $A$ passing through the center
$p$, not containing $\ell$. It is also not difficult to prove
that $\PEL(V)_p$ acts transitively on the lines passing through $p$
(here, one chooses an axis containing $p$, but with a different center).
The result follows from these observations.
\qed
\end{Lem}
\emph{Mutatis mutandis}, one proves that $\PEL(V)$ acts 2-transitively
on the hyperplanes.
\begin{Lem}
If $n\geq 2$, then $\EL(V)$ is perfect, i.e. $[\EL(V),\EL(V)]=\EL(V)$.
The same is true for $n=1$, provided that $|D|\geq 4$.

\proof
Assume first that $n\geq 2$. Given $\tau_{u\phi}$, choose $\psi$ linearly
independent from $\phi$, such that $\psi(u)=0$ and $v\in\kernel(\phi)$
with $\psi(v)=1$, Then $[\tau_{\psi u},\tau_{\phi v}]=\tau_{\phi u}$,
so every transvection is a commutator of transvections.
In the case $n=1$ one uses some clever matrix identities, and the fact
that the equation $x^2-1\neq 0$ has a solution in $D$ if $|D|\geq 4$,
see Hahn \& O'Meara \cite{HOM} 2.2.3.
\qed
\end{Lem}

\begin{Prop}
If $n\geq 2$ or if $n=1$ and $|D|\geq 4$, then $\PEL(V)$ is a simple
group.

\proof
The proof uses Iwasawa's simplicity criterion, see Hahn \& O'Meara
\cite{HOM} 2.2.B. The group $\PEL(V)$ is
perfect, and the stabilizer of a point $z$ has an abelian normal subgroup,
the group consisting of all elations with center $z$. These are the main
ingredients of the proof, see \emph{loc.cit.} 2.2.13.
See also Artin \cite{Artin} Ch.~IV Thm.~4.10.
\qed
\end{Prop}

We collect some further results about the group $\EL(V)$.
Let $H\cong\GL_1(D)$
denote the subgroup of $\GL(V)$ consisting of all matrices
of the form
\[
\begin{pmatrix}
a \\
 & 1 \\
 && 1 \\
 &&& 1 \\
 &&&& \ddots \\
 &&&&& 1
\end{pmatrix} 
\]
for $a\in D^\times$. One can show that $\GL(V)=H\EL(V)$, see
Hahn \& O'Meara \cite{HOM} 1.2.10 (the theorem applies, since
a skew field is a euclidean ring in the terminology of \emph{loc.cit.}).
An immediate consequence is the following.
\begin{Num}
If $D$ is commutative, then $\EL(V)=\SL(V)$.
\end{Num}
(Artin \cite{Artin} denotes the group $\EL(V)$ by $\SL(V)$, even if
$D$ is not commutative. Hahn \& O'Meara \cite{HOM} --- and other modern
books --- have a different
terminology; in their book, $\SL(V)$ is the kernel of the reduced norm.)

The next result follows from the fact that $\EL(V)$ acts strongly
transitively on the building $\Delta(V)$; such an action is always
primitive on the vertices of a fixed type, because the maximal parabolics
are maximal subgroups. In our case, the vertices of
the building $\Delta(V)$ are precisely the subspaces of $V$,
see Section \ref{AnBuilding}.
\begin{Num}
The action of $\EL(V)$ on $\Gr_k(V)$ is \emph{primitive}, for $1\leq k\leq n$.
It is \emph{two-transitive} if and only if $k=1,n$.
\end{Num}
Finally, we mention some exceptional phenomena.
\begin{Num}
Suppose that $D\cong\FF_q$ is finite, and let
$\PSL_m(q)=\PEL_m(\FF_q)$. There are the following isomorphisms
(and no others, see Hahn \& O'Meara \cite{HOM} p.~81).
\begin{gather*}
\PSL_2(2)\cong\Sym(3)\\
\PSL_2(3)\cong\Alt(4)\\
\PSL_2(4)\cong\PSL_2(5)\cong\Alt(5)\\
\PSL_2(7)\cong\PSL_3(2)\\
\PSL_2(9)\cong\Alt(6)\\
\PSL_4(2)\cong\Alt(8)
\end{gather*}
In particular, the groups $\PSL_2(q)$ are not perfect for $q=2,3$.
Note that the groups $\PSL_3(4)$ and $\PSL_4(2)$ have the same
order $20\,160$ without being isomorphic.
\end{Num}

\subsection{K$_1$ and the Dieudonn\'e determinant}
\begin{quote}
\emph{We explain the connection between the elementary linear group,
the first K-group $\sK_1(D)$ and the Dieudonn\'e determinant.}
\end{quote}
Recall the category $\Mod^\fin_D$  of all finite dimensional
vector spaces over $D$. Let
$\mathcal{M}=\{D^n|\ n\geq 0\}$; these vector spaces together with
the linear maps between them form a small and full subcategory of
$\Mod^\fin_D$.
Every vector space $V$ in $\Mod^\fin_D$ is isomorphic to a unique element
$[V]\in\mathcal{M}$. We make $\mathcal{M}$ into an additive semigroup
with addition $[V]+[W]=[V\oplus W]$, and neutral element $[0]$.
The dimension functor yields an isomorphism
\[
\xymatrix{{(\mathcal{M},+)}\ar[r]^\dim&{(\mathbb{N},+).}}
\]
Now $\sK_0(D)=\sK_0(\Mod^\fin_D)$ is defined to be the
Grothendieck group generated by
the additive semigroup of isomorphism classes of finite dimensional vector
spaces over $D$, i.e. $\sK_0(D)\cong\ZZ$.
(The Grothendieck group of a commutative semigroup $S$ is the universal
solution $G(S)$ of the problem
\[
\xymatrix{
S \ar[r]\ar[dr]_f & {G(S)} \ar @{..>}[d] \\
  & {H} }
\]
where $H$ is any group and $f$ is any semigroup homomorphism.)
While all this is rather trivial general nonsense
(but only since $D$ is a skew field!),
there are higher-rank K-groups which bear more information even for skew
fields.

For $n\leq m$, there is a natural inclusion
$\xymatrix@1{{\GL_n(D)\,}\ar@{^{(}->}[r]&{\GL_m(D)}}$
(as block matrices in the upper left, with $1$s on the diagonal
in the lower right).
Let $\GL_\stable(D)$ denote the direct limit over these inclusions
(this is called the \emph{stable linear group} --- not to be
confused with stability theory in the model theoretic sense), and let
$\EL_\stable(D)$ denote the corresponding direct limit over the groups
$\EL_n(D)$. For $n\geq 1$, there are exact sequences
\[
\xymatrix{
1 \ar[r] & {\EL_n(D)} \ar[r]\ar[d] & {\GL_n(D)} \ar[r]\ar[d] &
{\GL_n(D)/\EL_n(D)} \ar[r]\ar[d]^\cong & 1 \\
1 \ar[r] & {\EL_\stable(D)} \ar[r] & {\GL_\stable(D)} \ar[r] 
& {\GL_\stable(D)/\EL_\stable(D)} \ar[r] & 1 }
\]
We denote these quotients
\[
\sK_{1,n}(D)=\GL_n(D)/\EL_n(D)
\]
and put $\sK_1(D)=\GL_\stable(D)/\EL_\stable(D)$. The groups
$\sK_{1,n}(D)$ are
\emph{stable}, i.e. independent of $n$; there are isomorphisms
$\sK_{1,2}(D)\cong\sK_{1,3}(D)\cong\cdots\cong\sK_1(D)$, see
Hahn \& O'Meara \cite{HOM} 2.2.4. In this way we have
obtained the \emph{first K-group} $\sK_1(D)$.

\begin{Prop}
There is an isomorphism $\sK_1(D)\cong D^\times/[D^\times,D^\times]$. The
composite
\[
\xymatrix{
{\GL_n(D)} \ar[r]\ar[drr]^(0.6){\det} & {\GL_\stable(D)} \ar[d] \\
 & {\sK_1(D)} \ar[r]^(0.4)\cong & {D^\times/[D^\times,D^\times]} }
\]
is precisely the Dieudonn\'e determinant $\det$,
see Hahn \& O'Meara \cite{HOM}
2.2.2 and Artin \cite{Artin} Ch.~IV Thm.~4.6.
\qed
\end{Prop}
Since the determinant takes values in an abelian group, 
it is invariant under base change (matrix conjugation);
in particular, there is a well-defined (base independent)
determinant map
\[
\xymatrix{{\GL(V)}\ar[r]^(.4)\det&{D^\times/[D^\times,D^\times].}}
\]
We end this section with a commutative diagram which compares the
linear and the projective actions. Given $c\in D^\times$, we denote its
image in $D^\times/[D^\times,D^\times]$ by $\bar c$. Let
$C=\Cen(D^\times)$.
\[
\xymatrix{
& 1 \ar[d] & 1 \ar[d] & 1 \ar[d] \\
1 \ar[r] & {\{c\in C|\ \bar c^{n+1}=1\}} 
 \ar[r]\ar[d] & C   \ar[r]\ar[d]
& {\{\bar c^{n+1}|\ c\in C\}} \ar[r]\ar[d] & 1 \\
1 \ar[r] & {\EL(V)} \ar[r]\ar[d] & {\GL(V)} \ar[r]\ar[d] 
& {\sK_1(D)} \ar[r]\ar[d] & 1 \\
1 \ar[r] & {\PEL(V)} \ar[r]\ar[d] & {\PGL(V)} \ar[r]\ar[d] 
& {\sK_1(D)/\{\bar c^{n+1}|\ c\in C\}} \ar[r]\ar[d] & 1 \\
& 1 & 1 & 1 
}
\]
Here are some examples. If $D$ is commutative, then $\sK_1(D)=D^\times$.
Now let $\HH$ denote the quaternion division
algebra over a real closed field $\RR$. The \emph{norm} $N$ is defined as
$N(x_0+\textbf{i}x_1+\textbf{j}x_2+\textbf{ij}x_3)=
x_0^2+x_1^2+x_2^2+x_3^2$. Then $\xymatrix@1{{\HH^\times}\ar[r]^N&{\RR_{>0}}}$
has kernel $[\HH^\times,\HH^\times ]=\SS^3$, whence
$\sK_1(\HH)\cong\RR_{>0}$.

\subsection{Steinberg relations and K$_2$}\label{K2}
\begin{quote}
\emph{We introduce the Steinberg relations, which are the basic commutator
relations for projective spaces, and indicate briefly the connection with
higher K-theory.}
\end{quote}
In this section we assume that $\dim(V)=n+1$ is finite, and that $n\geq 2$.
Let $b_0,\ldots,b_n$ be a base for $V$, and let
$\beta_0,\ldots,\beta_n\in V^\vee$
be the dual base (i.e. $\beta_i(b_j)=\delta_{ij}$). Let
\[
\tau_{ij}(a)=\tau_{b_ia\beta_j},\qquad\text{ for }i\neq j.
\]
Thus $\tau_{ij}(a)$ can be pictured as a matrix with $1$s on the diagonal,
an entry $a$ at position $(i,j)$ ($i$th row, $j$th column), and $0$s else.
We claim that the elations $\{\tau_{ij}(a)|\ a\in D,\ i\neq j\}$
generate $\EL(V)$. Indeed, it is easy to see that the group generated
by these elations acts transitively on incident point-hyperplane pairs;
therefore, it contains all elations.
The maps $\tau_{ij}(a)$ satisfy the following relations, as is easily
checked.
\begin{description}
\item[SR1]
$\tau_{ij}(a)\tau_{ij}(b)=\tau_{ij}(a+b)$ for $i\neq j$.
\item[SR2]
$[\tau_{ij}(a),\tau_{kl}(b)]=1$ for $i\neq k$ and $j\neq l$.
\item[SR3]
$[\tau_{ij}(a),\tau_{jk}(b)]=\tau_{ik}(ab)$
for $i,j,k$ pairwise distinct.
\end{description}
These are the \emph{Steinberg relations}. 
They show that the algebraic structure of the skew field
$D$ is encoded in the little projective group $\PEL_n(D)$.

For each pair $(i,j)$ with $i\neq j$ we fix an isomorphic copy $U_{ij}$ of
the additive group $(D,+)$, and an isomorphism
$\xymatrix@1{{\tau_{ij}:(D,+)}\ar[r]^(.65)\cong&{U_{ij}}}$.
For $n\geq 2$, we define $\mathsf{St}_{n+1}(D)$ as the
free amalgamated product the
of $n(n+1)$ groups $U_{ij}=\{\tau_{ij}(a)|\ a\in D\}$,
factored by the normal subgroup generated by
the Steinberg relations \textbf{SR2}, \textbf{SR3}.
There is a natural epimorphism
\[
\xymatrix{{\mathsf{St}_{n+1}(D)}\ar[r]&{\EL_{n+1}(D)}}
\]
whose kernel is denoted $\sK_{2,n+1}(D)$. Again, there are natural
maps 
\[
\xymatrix{{\sK_{2,n+1}(D)}\ar[r]&{\sK_{2,m+1}(D)}}
\]
for $m\geq n$, and one can consider the
limits, the \emph{stable groups} $\mathsf{St}_\stable(D)$ and $\sK_2(D)$.
Clearly, there are exact sequences
\[
\xymatrix{
1\ar[r] & {\sK_{2,n+1}(D)} \ar[r] & {\mathsf{St}_{n+1}(D)} \ar[r] &
{\GL_{n+1}(D)} \ar[r] & {\sK_{1,n+1}(D)} \ar[r] & 1 }
\]
(and similarly in the limit).
The groups $\sK_2(D)$ bear some information about the skew field $D$.
See Milnor \cite{Milnor} for $\sK_2(D)$ of certain fields $D$;
for quaternion algebras, $\sK_2(D)$ is determined in Alperin-Dennis
\cite{AlpDen}.
One can prove that $\sK_{2,n+1}(D)\cong\sK_2(D)$ for $n>1$, see
Hahn \& O'Meara \cite{HOM} 4.2.18 and that
$\mathsf{St}_{n+1}(D)$ is a universal central extension of $\EL_{n+1}(D)$, 
provided that $n\geq 4$, see \emph{loc.cit.} 4.2.20,
or that $n\geq 3$ and that $\Cen(D)$ has
at least $5$ elements, see Strooker \cite{Strooker} Thm.~1.
We just mention the following facts.

(1) If $D$ is finite, then $\sK_2(D)=0$, so the groups
$\SL_{n+1}(q)$ are centrally closed for $n\geq 4$
(they don't admit non-trivial
central extensions), see Hahn \& O'Meara \cite{HOM} 2.3.10
(in low dimensions over small fields, there are exceptions, see
\emph{loc.cit.}). Also, the Steinberg relations yield a presentation of
the groups $\SL_{n+1}(q)$ for $n\geq 4$.

(2) Suppose that $D$ is a field with a primitive $m$th root of unity.
Let $\mathsf{Br}(D)$ denote its
Brauer group. Then there is an exact sequence of abelian groups
\[
\xymatrix{
1\ar[r] & {[\sK_2(D)]^m} \ar[r] & {\sK_2(D)} \ar[r] & {\mathsf{Br}(D)} \ar[r]
 & {[\mathsf{Br}(D)]^m} \ar[r] & 1}.
\]
(where we write $[A]^m=\{a^m|\ a\in A\}$ for an abelian group $(A,\cdot)$),
see Hahn \& O'Meara \cite{HOM} 2.3.12.

Applications of $\sK_2$, e.g.~in number theory, are mentioned in 
Milnor \cite{Milnor} and in Rosenberg \cite{Ros}.

\subsection{Different notions of projective space, characterizations}
\label{ProGeo}
\begin{quote}
\emph{We introduce the point-line geometry and the building obtained from
a projective geometry and compare the resulting structures.
Then we mention the 'second' Fundamental Theorem of Projective Geometry
which characterizes projective geometries of rank at least $3$.
We describe the Tits system (BN-pair) for the projective geometry and
the root system, and we explain how the root system reflects properties
of  commutators of root elations.}
\end{quote}
\label{AnBuilding}
A \emph{point-line geometry} is a structure
\[
(\cP,\cL,*),
\]
where $\cP$ and
$\cL$ are non-empty disjoint sets, and
$*\subseteq (\cP\cup\cL)\times (\cP\cup\cL)$
is a symmetric and reflexive
binary relation, such that $*|_{\cP\times\cP}=\id_\cP$ and
$*|_{\cL\times\cL}=\id_\cL$.
Given a projective geometry $\PG(V)=(\Gr_1(V),\ldots,\Gr_n(V),*)$ of rank
$n\geq 2$, we can consider the point-line geometry
\[
\PG(V)_{1,2}=(\Gr_1(V),\Gr_2(V),*).
\]
It is easy to recover the whole structure $\PG(V)$ from this;
call a set $X$ of points a \emph{subspace} if it has the following
property: for every triple of pairwise distinct collinear points $p,q,r$
(i.e. there exists a line $\ell\in\Gr_2(V)$ with $p,q,r *\ell$), we have
the implication
\[
(p,q\in X)\quad\Longrightarrow\quad(r\in X).
\]
We define the \emph{rank} of a subspace inductively as $\rk(\emptyset)=-1$,
and $\rk(X)\geq k+1$ if $X$ contains a proper subspace $Y\subsetneq X$
with $\rk(Y)\geq k$. Clearly, $\Gr_{k+1}(V)$ can be identified with the
set of all subspaces of rank $k$.

The point-line geometry $(\cP,\cL,*)=(\Gr_1(V),\Gr_2(V),*)$ has the following
properties.
\begin{description}
\item[PG1]
Every line is incident with at least $3$ distinct points.
\item[PG2]
Any two distinct points $p,q$ can be joined by a unique line which
we denote by $p\vee q$.
\item[PG3]
There exist at least $2$ distinct lines.
\item[PG4]
If $p,q,r$ are three distinct points, and if $\ell$ is a line which
meets $p\vee q$ and $p\vee r$ in two distinct points,
then $\ell$ meets $q\vee r$.
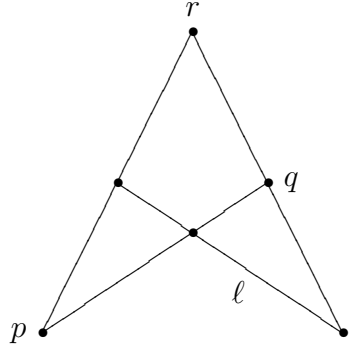
\[
\begin{xy} <2cm,0cm>:
(0,0)*@{*}="a"*++!R{p},
(2,0)*@{*}="b",
(1,2)*@{*}="c"*++!D{r},
(0.5,1)*@{*}="d",
(1.5,1)*@{*}="e"*++!L{q},
"a";"c"**@{-},
"b";"c"**@{-},
"b";"d"**@{-}?(0.4)*+!UR{\ell},
"a";"e"**@{-} ?!{"b";"d"} *@{*}="f"
\end{xy}
\]
\end{description}
Axiom \textbf{PG4} is also called \emph{Veblen's axiom} or
the \emph{Veblen-Young property}.
A point-line geometry which satisfies these axioms is called a
\emph{projective (point-line) geometry}. If there exist two lines
which don't intersect,
then it is called a \emph{projective space}, otherwise a
\emph{projective plane}.

Veblen's axiom PG4 is the important 'geometric' axiom in this list;
the axioms PG1--PG3 exclude only some obvious pathologies. It is one
of the marvels of incidence geometry that this simple axiom encodes 
--- by the Fundamental Theorem of Projective Geometry stated below ---
the whole theory of skew fields, vector spaces, and linear algebra.

\begin{Num}\textbf{The Fundamental Theorem of Projective Geometry, II}
\label{FTPG2}\

\noindent\em
Let $(\cP,\cL,*)$ be a projective (point-line) space which is not a
projective plane. Then there exists a
skew field $D$, unique up to isomorphism, and a right vector
space $V$ over $D$, unique up to isomorphism, such that 
\[
(\cP,\cL,*)=(\Gr_1(V),\Gr_2(V),*).
\]
Here, the vector space dimension can be infinite. In fact, the dimension
is finite if and only if one the following holds:

(1) Every subspace has finite rank.

(2) There exists no subspace $U$ and no automorphism $\phi$ such that
$\phi(U)$ is a proper subset of $U$.

\medskip\noindent\em
The theorem is folklore; we just refer to L\"uneburg \cite{Lun1}
\cite{Lun2}, or to Faure \& Fr\"olicher \cite{FauFro1}
for a category-theoretic proof.
\qed
\end{Num}
Let $V^\vee$ denote the dual of $V$. This is a right vector space
over $D^\op$; it is easy to see that there is an isomorphism
\[
\PG(V^\vee)_{1,2}\cong(\Gr_n(V),\Gr_{n-1}(V),*).
\]
Now consider the following structure (for finite dimension $\dim(V)=n+1$).
Let $\mathcal{V}=\Gr_1(V)\cup\Gr_2(V)\cup\cdots\cup\Gr_n(V)$ and let
$\Delta$ denote the collection of all subsets of $\mathcal{V}$ which
consist of pairwise incident elements. Such a set is finite and
has at most $n$ elements. Then $\Delta(V)=(\Delta,\subseteq)$
is a poset, and in fact an abstract simplicial complex of dimension
$n-1$.
The set $\mathcal{V}$ can be identified with the minimal elements
(the vertices) of $\Delta$. There is an exact sequence
\[
\xymatrix{
1 \ar[r]  & {\Aut(\PG(V))} \ar[r] & {\Aut(\Delta(V))} \ar[r]
 & {\mathsf{AAut}(D)/\Aut(D)} \ar[r] & 1 }
\]
The poset $\Delta(V)$ is the \emph{building} associated to the
projective
space $\PG(V)$. 
See Brown \cite{Brown}, Garrett \cite{Garrett},
Grundh\"ofer \cite{Theo} (these proceedings),
Ronan \cite{Ronan}, Taylor \cite{Taylor}, Tits \cite{TitsLNM}.
Here, we view a building as a simplicial complex, without a
type function (the type function $\mathsf{type}$
would associate to a vertex
$v\in\Gr_k(V)$ the number $k$). 
Such a type function can always be defined
and is unique up to automorphisms; if we consider only type-preserving
automorphisms, then we obtain $\Aut(\PG(V))$ as the automorphism group.
(In Grundh\"ofer \cite{Theo}, the buildings are always endowed with a
type funtion.)
There is no natural way to recover $\PG(V)$ from
$\Delta(V)$, but we can recover both $\PG(V)$ and $\PG(V^\vee)$
simultaneously. The following diagram shows the various 'expansions' and
'reductions'. Only the solid arrows describe natural constructions; the
dotted arrows require the choice of a type function, i.e. one has
to choose which elements are called points, and which ones
are called hyperplanes.
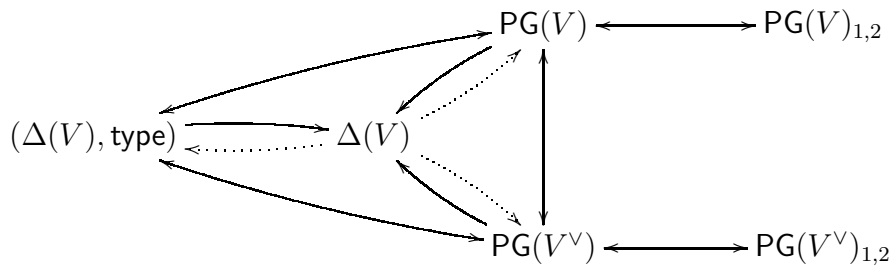
\[
\xymatrix{
&&& {\PG(V)} \ar @/_1ex/[dl]\ar @{<->}[rr]\ar @{<->}[dd]
&& {\PG(V)_{1,2}}  \\
{(\Delta(V),\mathsf{type})} \ar @{<->}@/^1ex/[urrr] \ar @{<->}@/_1ex/[drrr]
\ar @/^1ex/[rr]
&& {\Delta(V)} \ar @{..>}@/_1ex/[ur] \ar @{..>}@/^1ex/[dr] 
\ar @{..>}@/^1ex/[ll]\\
&&& {\PG(V^\vee)} \ar @/^1ex/[ul] \ar @{<->}[rr]
&& {\PG(V^\vee)_{1,2}} 
}
\]
This is probably the right place to introduce the \emph{Tits system}
(or \emph{BN-pair})
of $\PEL(V)$. Actually, it is easier if we lift everything into
the group $\EL(V)$ (we could equally well work with the Steinberg
group $\mathsf{St}_n(D)$, or the general linear group $\GL(V)$).
\begin{Num}\textbf{The Tits system for $\EL(V)$}
Let $b_0,\ldots,b_n$ be a base for $V$, and let $p_i=b_iD$.
Every proper subset of $\{p_1,\ldots,p_n\}$ spans a subspace; in this
way, we obtain a collection 
\[
\Sigma^{(0)}=\left.\left\{V_J=\mathsf{span}\{b_j|\ j\in J\}\right| 
\emptyset\neq J\subsetneq \{0,\ldots,n\}\right\}
\]
of $2(2^n-1)$ subspaces.
With the natural inclusion '$\subseteq$', this becomes a poset
and an abstract simplicial complex; combinatorically, this is the
complex $\partial\Delta^{n+1}$
of all proper faces of a $n+1$-simplex; for $n=2$, we have the
points and sides of a triangle, and for $n=3$ the points, edges and
sides of a tetrahedron.

Now we consider a different simplicial complex, the \emph{apartment}
$\Sigma$. The vertices of $\Sigma$ are the elements of $\Sigma^{(0)}$,
and the higher rank elements are sets of pairwise incident elements.
This complex can be pictured as follows: consider the first barycentric
subdivision $\mathsf{Sd}\Sigma^{(0)}$ of $\Sigma^{(0)}$; the barycentric
subdivision adds a vertex in every face of $\Sigma^{(0)}$.
This flag complex 
\[
\Sigma=\mathsf{Sd}\partial\Delta^{n+1}
\]
is the \emph{apartment} spanned by $p_0,\ldots,p_n$.

Let $T\subseteq\EL(V)$ denote the pointwise stabilizer of $p_0,\ldots,p_n$
(equivalently, the elementwise stabilizer of $\Sigma^{(0)}$ or $\Sigma$),
and $N$ the setwise stabilizer of $\{p_0,\ldots,p_n\}$
(or $\Sigma^{(0)}$, or $\Sigma$). Thus
$N/T\cong\Sym(n+1)$; this quotient is the \emph{Weyl group} for $\PGL(V)$.
Let $B\subseteq\EL(V)$ denote the stabilizer of the flag
\[
C=(p_0,p_0\oplus
p_2,p_0\oplus p_1\oplus p_2,\ldots,p_0\oplus\cdots\oplus p_{n-1})
\]
With respect to the base $b_0,\ldots,b_n$, the group
$B$ consists of the upper
triangular matrices in $\EL(V)$, and $N$ consists of permutation 
(or monomial) matrices
(a permutation matrix has precisely one non-zero entry in every row and
column), and $T=B\cap N$ consists of diagonal matrices with the property
that the product of the entries lies in the commutator group
$[D^\times,D^\times]$.
For $1\leq i\leq n$, let $s_i$ be the $T$-coset of the matrix
\[
\begin{pmatrix}
1\\
& \ddots \\
 && 1 \\
 &&& 0 & 1 \\
 &&& -1& 0 \\
 &&&&& 1 \\
 &&&&&& \ddots \\
 &&&&&&& 1 
\end{pmatrix} 
\]
which interchanges $p_i$ and $p_{i-1}$. Then $s_i$ is an involution in $W$,
and 
\[
(W,\{s_1,\ldots,s_n\})
\]
presents $W$ as a \emph{Coxeter group}.
It is a routine matter to check that these data
$(G,B,N,\{s_1,\ldots,s_n\})$ satisfy the axioms of a \emph{Tits system}:
\begin{description}
\item[TS1]
$B$ and $N$ generate $G$.
\item[TS2]
$T=B\cap N$ is normalized by $N$.
\item[TS3]
The set $S=\{s_1,\ldots,s_n\}$ generates $N/T$ and has the following
properties.
\item[TS4]
$sBs\neq s$ for all $s\in S$.
\item[TS5]
$BsBwB\subseteq BwB\cup BswB$ for all $s\in S$ and $w\in W=N/T$.
\end{description}
The group $T\subseteq B$ has a normal complement, the group 
$U$ generated by all elations $\tau_{ij}(a)$, for $i<j$.
Thus $B$ is a semidirect product $B=TU$, and $(\EL(V),B,N,S)$ is what
is called a (strongly split) \emph{Tits system} (or \emph{BN-pair}).
(A Tits system is called \emph{split} if $B$ can be written as a
not necessarily semidirect product $B=TU$, such that $U\lhd B$ is
normal in $B$. Usually, one also requires $U$ to be nilpotent ---
this is the case in our example.
The group $U$ acts transitively on the apartments
containing the chamber corresponding to $B$. If $U\cap T=1$, the
Tits system is said to be \emph{strongly split}.)
Note also that $B$ is \emph{not} solvable if $D$ is not commutative.
\end{Num}
Now we describe the root system for $\Delta(V)$.
\begin{Num}\textbf{Thin projective spaces and the root system}
Let $\cP=\{0,\ldots,n\}$ and let $\cL=\binom{\cP}2$ denote the collection
of all $2$-element subsets of $\cP$. The incidence $*$ is the
symmetrized inclusion relation. Then $(\cP,\cL,*)$ is a
thin projective geometry of rank $n$, i.e. a projective space where
every line is incident with precisely two points. The corresponding
thin projective space (the projective geometry over the 'field $0$
with one element') is
\[
\textstyle
\PG_n(0)=\left(\binom\cP1,\binom\cP2,\ldots,\binom\cP n,\subseteq\right).
\]
Now we construct a 'linear model' for this geometry.
Consider the real euclidean vector space $\RR^{n+1}$ with its standard
inner product $\bra{-,-}$ and standard base $e_0,e_1,\ldots,e_n$.
Let $E$ denote the orthogonal complement
of the vector $v=e_0+\cdots+e_n=(1,1,\ldots,1)\in\RR^{n+1}$.
Let 
\[
\textstyle p_i=e_i-v\frac1n\in E,\text{ for }i=0,\ldots,n.
\]
We identify $p_0,\ldots,p_n$ with the points of $\PG_n(0)$.
The subspaces of rank $k$ correspond precisely to (linearly independent)
subsets $p_{i_1},\ldots,p_{i_k}$ of $\cP$.
We identify such a subspace with the vector $\frac1k(p_{i_1}+\cdots+p_{i_k})$;
this is the barycenter of the convex hull of $\{p_{i_1},\ldots,p_{i_k}\}$.
This is our model of $\PG_n(0)$; the picture shows the case $n=2$.
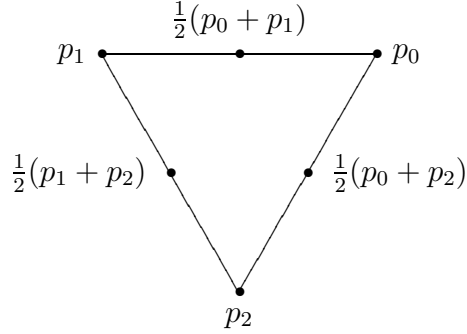
\[
\begin{xy} /r5pc/:
{\xypolygon3"R"{~={30}~<{}~><{@{-}}~>>{|*@{*}}}},
"R1"*++!L{p_0},
"R3"*++!U{p_2},
"R2"*++!R{p_1},
"R1"*@{*},
"R3"*@{*},
"R2"*@{*},
"R1";"R2"**{}?(.5)*++!D{\frac12(p_0+p_1)},
"R1";"R3"**{}?(.5)*+++!L{\frac12(p_0+p_2)},
"R2";"R3"**{}?(.5)*+++!R{\frac12(p_1+p_2)}
\end{xy}
\]
Now we construct the Weyl group.
For $i\neq j$ put $\eps_{ij}=e_i-e_j$. The reflection $r_{ij}$
at the hyperplane $\eps_{ij}^\perp$ in $E$,
\[
\xymatrix{x\ar@{|->}[r]^(.3){r_{ij}}&{x-\bra{x,\eps_{ij}}\eps_{ij}},}
\]
permutes $p_0,\ldots,p_n$; the isometry group $W$ generated by these
reflections is the \emph{Weyl group} of type $\mathbf A_n$
(which is isomorphic to the Coxeter group $N/T$ of the Tits system).

These vectors $\eps_{ij}$
form a \emph{root system of type $\mathbf A_n$} in $E$, with
$\Phi=\{\eps_{ij}|\ i\neq j\}$ as set of \emph{roots}. We call
the set $\Phi^+=\{\eps_{ij}|\ i<j\}$ the set of \emph{positive roots}; this
determines the $n$ \emph{fundamental roots}
$\{\eps_{01},\eps_{12},\ldots,\eps_{n-1,n}\}$. The fundamental roots form a base
of $E$, and every root is an integral linear combination of fundamental roots,
such that either all coefficients are non-negative (this yields the positive
roots) or non-positive.

To each root $\eps_{ij}$, we attach the group $U_{ij}$ as defined in 
Section \ref{K2}.
From the Steinberg relations, we see the following: if $i<j$ and $k<l$, then
$[U_{ij},U_{kl}]=0$ if there exists no positive root
which is a linear combination $\eps_{ij}a+\eps_{kl}b$, with $a,b\in\ZZ_{>0}$.
For our root system, the only instance where such a linear combination is
a positive root is when $\mathsf{card}\{i,j,k,l\}=3$, and in
this case the Steinberg relations show that the commutator is indeed
not trivial. The picture below shows the root system $\mathbf A_2$
(the case $n=2$); the fundamental roots are $\eps_{01}$ and $\eps_{12}$, and
the positive roots are $\eps_{01}$, $\eps_{12}$,
and $\eps_{02}=\eps_{01}+\eps_{12}$.
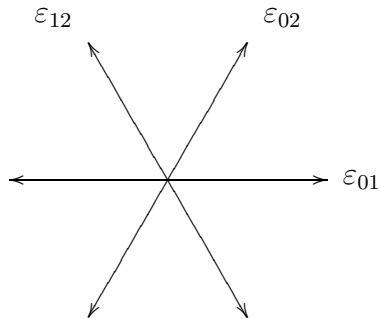
\[
\begin{xy} /r5pc/:
{\xypolygon6"R"{~<<{@{->}}~>{}}},
"R1"*++!L{\eps_{01}},
"R2"*++!DL{\eps_{02}},
"R3"*++!DR{\eps_{12}}
\end{xy}
\]
The hyperplanes $\eps_{ij}^\perp\subseteq E$ yield a triangulation of
the unit sphere $\SS^{n-1}\subseteq E$; as a simplicial complex, this
is precisely the apartment $\Sigma$. The half-apartments correspond to
the half-spaces $\{v\in E|\ \bra{v,\eps_{ij}}\geq 0\}$. From this, it
is not hard to see that the groups $U_{ij}$ are root groups in the
building-theoretic
sense (as defined in Grundh\"ofer's article \cite{Theo}
in these proceedings), see also Ronan \cite{Ronan} Ch.~6.
\end{Num}

\subsection{The little projective group as a 2-transitive group}
\begin{quote}
\emph{We show that the projective space is determined by
(and can be recovered from) the action of the elementary linear group
on the point set.}
\end{quote}
Let $\PEL(V)\subseteq H\subseteq\PGL(V)$ be a subgroup, and assume that
$\dim(V)\geq 3$. Then $(H,\Gr_1(V))$ is a 2-transitive permutation group.
Let $L\subseteq\Gr_1(V)$ be a point row, i.e. the set of all points lying
on a line $\ell\in\Gr_2(V)$. Let $p,q\in L$ be distinct points. Then
$H_{p,q}$ has precisely four orbits in $\Gr_1(V)$: the two singletons
$\{p\},\{q\}$, the set $L\setminus\{p,q\}$, and $X=\Gr_1(V)\setminus L$.
The set $X$ has the property that every $h$ in $H$ which fixes $X$ pointwise
fixes $\Gr_1(V)$ pointwise. None of the other three orbits has this property.
Thus one can see the line $L\subset\Gr_1(V)$ from the $H$-action, we have
a canonical (re)construction
\[
\xymatrix{{(H,\Gr_1(V))}\ar@{|->}[r]&{(H,\Gr_1(V),\Gr_2(V),*).}}
\]
Combining this with the Fundamental Theorem of Projective Geometry
\ref{FTPG1}, we have
the next result.
\begin{Prop}\label{PermChar1}
For $i=1,2$, let $\PG(V_i)$ be projective geometries (of possibly
different ranks $n_i\geq 2$)
over skew fields $D_1,D_2$. Let $\PEL(V_i)\subseteq H_i\subseteq\PGL(V_i)$
be subgroups. If there exists an isomorphism of permutation groups
\[
\xymatrix{{(H_1,\Gr_1(V_1))}\ar[r]^\phi_\cong&{(H_2,\Gr_1(V_2))}}
\]
then there exists a semilinear bijection $\xymatrix@1{{F:V_1}\ar[r]&{V_2}}$
which induces $\phi$ (and $n_1=n_2$).
\qed
\end{Prop}
The result is also true in dimension $2$, but the proof is more complicated,
as we will see in Section \ref{ProjectiveLine}.
The problem whether an abstract group isomorphism
$\xymatrix@1{{H_1}\ar[r]_\cong^\phi&{H_2}}$
is always induced by a semilinear map is much more subtle. The result
is indeed that such an isomorphism is induced by a linear map, composed
with an isomorphism or anti-isomorphism of the skew fields in question,
provided that the vector space dimensions are large enough (at least $3$),
see Hahn \& O'Meara \cite{HOM} 2.2D.
The crucial (and difficult) step is to show that $\phi$ maps transvections
to transvections.

\subsection{Projective planes}
\begin{quote}
\emph{We mention the classification of Moufang planes.}
\end{quote}
The Fundamental Theorem of Projective Geometry does not apply to
projective planes. However, there is the following result.
Suppose that $(\cP,\cL,*)$ is a projective plane. Given a flag $(p,\ell)$,
(i.e. $p*\ell)$), the group $G_{[p,\ell]}$ is defined to be the set of 
all automorphisms which fix $p$ linewise and $\ell$ pointwise
(so for $\PG(D^3)$, we have
$G_{[p,\ell]}=U_{(p,\ell)}$ in our previous notation).
Let $h$ be a line passing through $p$ and different from $\ell$.
It is easy to see that $G_{[p,\ell]}$ acts freely on the set
$H'=\{q\in\cP|\ q\neq p,\ q*\ell\}$. If this action is transitive,
then $(\cP,\cL,*)$ is called \emph{$(p,\ell)$-homogeneous.}
The projective plane is called a \emph{Moufang plane} if it is
$(p,\ell)$-homogeneous for any flag $(p,\ell)$.
If $(\cP,\cL,*)=\PGL(V)$ for some $3$-dimensional vector space $V$,
then we have a Moufang plane.

Recall that an \emph{alternative field} is a (not necessarily associative)
algebra with unit, satisfying the following relations.
\begin{description}
\item[AF1]
If $a\neq 0$, then the equations $ax=b$ and $ya=b$ have unique solutions
$x,y$.
\item[AF2]
The equalities $x^2y=x(xy)$ and $yx^2=(yx)x$ hold for all $x,y$.
\end{description}
Clearly, every field or skew field is an alternative field.
The structure theorem of alternative fields says that every
non-associative alternative field is a central 8-dimensional
algebra over a field $K$, a so-called \emph{Cayley division algebra}.
Not every field $K$ admits a Cayley division algebra; it is necessary
that $K$ admits an anisotropic quadratic form of dimension $8$, so
finite fields or algebraically closed fields do not admit Cayley
division algebras. Every real closed (or ordered)
field admits a Cayley division algebra.

Given an alternative field $A$, we construct a projective plane
$\PG_2(A)$ as follows. Let $\infty$ be a symbol which is not
an element of $A$. Let
\begin{align*}
\cP&=\{(\infty)\}\cup\{(a)|\ a\in D\}\cup\{(x,y)|\ x,y\in A\} \\
\cL&=\{[\infty]\}\cup\{[a]|\ a\in D\}\cup\{[x,y]|\ x,y\in A\} \\
\end{align*}
The incidence $*$ is defined as
\[
(\infty)*[a]*(a,sa+t)*[s,t]*(s)*[\infty]*(\infty)
\]
Note that for a field or skew field $D$, this is precisely
$\PG(D^3)$.
\begin{Thm}[Moufang Planes]
\

\noindent
Let $(\cP,\cL,*)$ be a Moufang plane. The there exists an alternative
field $A$, unique up to isomorphism, such that
$(\cP,\cL,*)$ is isomorphic to $\PG_2(A)$.

\medskip\noindent\em
For a proof see Hughes \& Piper \cite{HP} or Pickert \cite{Pi};
the structure theorem for non-associative alternative fields is
proved in Van Maldeghem \cite{HVM}.
\qed
\end{Thm}
The Moufang planes can also be described by means of Steinberg
relations: for each pair $(i,j)$ with $i\neq j$ and $i,j\in\{1,2,3\}$, fix
a copy $U_{ij}\cong (A,+)$. Let $G$ denote the free product of these
six groups, factored by the Steinberg relations as given in \ref{K2}
(note that the Steinberg relations make sense even in the non-associative
case).
The group $G$ has a natural
Tits system which yields the Cayley plane $\PG_2(A)$. The group
induced by $G$ on $\PG_2(A)$ is a $K$-form of a simple adjoint algebraic
group of type $\mathbf E_6$.

\subsection{The projective line and Moufang sets}
\label{ProjectiveLine}
\begin{quote}
\emph{We investigate the action of the linear group as a permutation group
on the projective line. This is a special case of a Moufang set.}
\end{quote}
The projective line $\cP=\Gr_1(V)$, for $V\cong D^2$,
is a set without further structure. We add structure by
specifying properties of the group $G=\PGL(V)$ acting on it.
This group has two remarkable properties: 
(1) $G$ acts 2-transitively on $\cP$.
(2) The stabilizer $B=G_p$ of a point $p=vD\in\cP$ has a regular
normal subgroup $U_p$, the group induced by maps of the form
$\id_V+v\rho$, where $\rho$ runs through the collection of all
non-zero elements of $V^\vee$ which annihilate $v$.

Moufang sets where first defined by Tits in \cite{TitsDurham};
our definition given below is stated in a slightly different way.
We define a \emph{Moufang set} as a triple
$(G,U,X)$, where $G$ is a group acting on a set $X$ (with at least 3
elements), and $U$ is a subgroup of $G$.
We require the following properties.
\begin{description}
\item[MS1]
The action of $G$ on $X$ is 2-transitive.
\item[MS2]
The group $U$ fixes a point $x$ and acts regularly on $X\setminus\{x\}$.
\item[MS3]
The group $U$ is normal in the stabilizer $G_x$.
\end{description}
The properties \textbf{MS2} and \textbf{MS3} will be summarized in the
sequel as \emph{'$G_x$ has a regular normal subgroup'}; we will also
say that \emph{'$U$ makes $(G,X)$ into a Moufang set'}.
Let $y\in X\setminus\{x\}$, and put $T=G_{x,y}$.
Then clearly $G_x=TU$ is a semidirect product. (The pair
$(G_x,T)$ is a what is called a (strongly) split Tits system (BN-pair)
of rank 1 for the group $G$.)

If $U=G_x$, then $T=1$, so $G$ is sharply 2-transitive. This case
has its own, special flavor.
Note that in general, $U$ is \emph{not}
determined by $G$ and $x$. As a counterexample,
let $\HH$ denote the quaternion division algebra over a real closed
field, let $X=\HH$, and consider the group consisting of maps of the form
$[\xymatrix@1{x\ar@{|->}[r]&axb+t}]$, for $a,b\in \HH^\times$ and $t\in\HH$.
The stabilizer of $0$ 
consists of the maps
$[\xymatrix@1{x\ar@{|->}[r]&axb}]$, and it has two regular normal
subgroups
isomorphic to $\HH^\times$,
consisting of the maps $[\xymatrix@1{x\ar@{|->}[r]&ax}]$ or
$[\xymatrix@1{x\ar@{|->}[r]&xb}]$.

Let $H\subseteq\PGL_2(D)$ be a subgroup containing $\PEL_2(D)$.
We identify the projective line $\Gr_1(D^2)$ with $D\cup\{\infty\}$,
identifying $\binom{x}1D$ with $x$ and $\binom10D$
with $\infty$. Let $U_\infty$
denote the group consisting of the maps
$[\xymatrix@1{x\ar@{|->}[r]&x+t}]$, for $t\in D$.
So 
\[
(H,U_\infty,D\cup\{\infty\})
\]
is a Moufang set. The stabilizer
$T=H_{0,\infty}$ contains all maps
$[\xymatrix@1{x\ar@{|->}[r]&axa}]$, for $a\in D^\times$.
In particularly, $T$ is commutative if and only if $D$ is commutative.
Now we consider the following problem: 
\begin{center}
\emph{Is it possible to recover $U_\infty$
from the action of $H$?}
\end{center}
In the commutative case, the answer is easy:
$U_\infty$ is the commutator group of $H_\infty$,
\[
T\text{ commutative }\quad\Longrightarrow\quad
U_\infty=[H_\infty,H_\infty].
\]
Also, an element $1\neq g\in H_\infty$ is
contained in $U_\infty$ if and only if
$g$ has no fixed point in $D$. Thus, $U_\infty$ is the \emph{only} regular
normal subgroup of $H_\infty$.

Now suppose that $D$ is not commutative. In this case, we will prove first
that the action of $H_\infty$ on $D$ is primitive.
We have to show that $T$ is a maximal subgroup of $H_\infty$. 
Let $g\in H_\infty\setminus T$, and consider the group $K$ generated
by $T$ and $g$. Since $H_\infty$ splits as a semidirect product
$H_\infty=TU_\infty$,
we may assume that $g=[\xymatrix@1{x\ar@{|->}[r]&x+t}]\in U$, with $t\neq 1$.
Since $K$ contains $T$, it contains all maps of the form
$[\xymatrix@1{x\ar@{|->}[r]&{x+taba^{-1}b^{-1}}}]$.
Using some algebraic identities
for multiplicative commutators as in Cohn \cite{Cohn} Sec.~3.9, one
shows that
the set of all multiplicative commutators generates $D$ additively.
Thus $U_\infty\subseteq K$, whence $K=H_\infty$,
\begin{align*}
T\text{ not commutative } &\quad\Longrightarrow\quad
H_\infty\text{ primitive }\\
&\quad\Longrightarrow\quad
U_\infty\text{ unique abelian normal subgroup of }H_\infty
\end{align*}
(for the last implication see Robinson \cite{Rob} 7.2.6).
As in the commutative case,
there is no other way of making the projective line into a Moufang set.
Indeed, suppose that $U_\infty\neq U'\unlhd H_\infty$
is another regular normal subgroup.
Then $U_\infty\cap U'\unlhd H_\infty$ is also normal, so either
$U_\infty\cap U'=1$, or
$U'\supseteq U_\infty$.
In the latter case, $U_\infty=U'$ since we assumed the action to
be regular. So suppose $U_\infty\cap U'=1$. Then $U'U_\infty$ is a
direct product.
Define a map $\xymatrix@1{\phi:U_\infty\ar[r] & U'}$
by putting $\phi(u)=u'$ if and only
if
$u(0)=u'(0)$.
Then $\phi(u_1u_2)(0)=(u_1u_2)'(0)=u_1u_2'(0)=u_2'u_1(0)
=u_2'u_1'(0)$, so $\phi$ is an anti-isomorphism. In particular,
$U'$ is abelian, whence $U_\infty=U'$, a contradiction.
\begin{Prop}
Let $V$ be a two-dimensional vector space over a skew field $D$, and
assume that $\PEL(V)\subseteq H\subseteq \PGL(V)$. Then $G_x$ contains
a unique regular normal subgroup $U_\infty$, i.e. there is a unique
way of making $(H,\Gr_1(V))$ into a Moufang set.
\end{Prop}
We combine this with Hua's Theorem.
\begin{Thm}[Hua]
Let 
\[
\xymatrix{{(\PEL(V),U_\infty,D\cup\{\infty\})}\ar[r]^\alpha_\cong
&{(\PEL(V'),U_\infty',D'\cup\{\infty'\})}}
\]
be an isomorphism of Moufang sets. Then $\alpha$ is induced by an
isomorphism or anti-isomorphism of skew fields.

\medskip\noindent\em
For a proof see Tits \cite{TitsLNM} 8.12.3 or Van Maldeghem \cite{HVM}
p.~383--385.
\end{Thm}
\begin{Cor}
Let $V,V'$ be 2-dimensional vector spaces over skew fields $D,D'$, let
$\PEL(V)\subseteq H\subseteq\PGL(V)$ and 
$\PEL(V')\subseteq H'\subseteq\PGL(V')$ and assume that
\[
\xymatrix{{(H,\Gr_1(V))}\ar[r]^\alpha_\cong&{(H',\Gr_1(V'))}}
\]
is an isomorphism of permutation groups. Then there exists an isomorphism
or anti-isomorphism
$\xymatrix@1{D\ar[r]^\theta_\cong&{D'}}$ which induces $\alpha$.

\proof
Clearly, $\alpha$ induces an isomorphism of the commutator groups
$[H,H]=\PEL(V)$ and $[H',H']=\PEL(V')$ (we can safely disregard the
small fields $\FF_2$ and $\FF_3$, since here, counting suffices).
There is a unique way of making these permutation groups into Moufang
sets, and to these Moufang sets, we apply Hua's Theorem.
\qed
\end{Cor}
The following observation is due to Hendrik Van Maldeghem.
Let $V$ be a vector space of dimension at least $3$, let $H$ be a group
of automorphisms of $\PG(V)$ containing $\PEL(V)$. Then $(H,\Gr_1(V))$
cannot be made into a Moufang set.
To see this, let $p\in\Gr_1(V)$ and assume that $U\unlhd H_p$ is
a normal subgroup acting regularly on $\Gr_1(V)\setminus\{p\}$.
Let $u\in U$, and let $\tau$ be an elation with center $p$. Then
$u\tau u^{-1}$ is also an elation with center $p$, and so is
$[u,\tau]$. If we choose $u$, $\tau$ in such a way that $u$ doesn't fix
the axis of $\tau$ (which is possible, since $\dim(V)\geq 3$), then
$[u,\tau]\in U$ is a non-trivial elation with center $p$. Since $U$ is normal,
$U$ contains all elations with center $p$. These elations form an abelian
normal subgroup of $H_p$ which is, however, not regular on $\Gr_1(V)$.
\begin{Lem}
Let $\PEL(V)\subseteq H\subseteq \PGL(V)$ and assume that $\dim(V)\geq 3$.
Then $H_p$ contains no regular normal subgroup; in particular,
$(H,\Gr_1(V))$ cannot be made into a Moufang set.
\qed
\end{Lem}
Combining the results of this section with Proposition \ref{PermChar1},
we have the following
final result about actions on the point set.
\begin{Cor}
For $i=1,2$, let $V_i$ be vector spaces over skew fields $D_i$, of
(finite)
dimensions $n_i\geq 2$. Let $\PEL(V_i)\subseteq H_i\subseteq\PGL(V_i)$
and assume that
\[
\xymatrix{{(H_1,\Gr_1(V_1))}\ar[r]^\phi_\cong&{(H_2,\Gr_1(V_2))}}
\]
is an isomorphism of permutation groups. Then $n_1=n_2$, and
$\phi$ is induced by a semilinear isomorphism, except if
$n_1=2$, in which case $\phi$ may also be induced by an anti-isomorphism
of skew fields.
\end{Cor}
These results are  also true for infinite dimensional vector spaces.

\section{Polar spaces and quadratic forms}
In this second part we consider $(\sigma,\eps)$-hermitian forms
and their generalizations, pseudo-quadratic forms. From now on,
we consider
also infinite dimensional vector spaces. As a finite dimensional
motivation, we start with dualities.
Suppose that $\dim(V)=n+1$ is finite, and that
there is an isomorphism $\phi$ between the projective
space $\PG(V)$ and its dual,
\[
\xymatrix{{\PG(V)}\ar[r]^\phi_\cong&{\PG(V^\vee).}}
\]
Now $V^\vee$ is in a natural way
a right vector space over the opposite skew field
$D^\op$.
By the Fundamental Theorem of Projective Geometry \ref{FTPG1},
$\phi$ is
induced by a
$\sigma$-semilinear bijection $\xymatrix@1{{f:V}\ar[r]^(.6)\cong&{V^\vee}}$,
relative to an isomorphism $\xymatrix@1{{\sigma:D}\ar[r]&{D^\op}}$, i.e.
$\sigma$ is an
anti-automorphism of $D$.
Note also that there is a natural isomorphism
\[
\PG(V^\vee)\cong(\Gr_n(V),\Gr_{n-1}(V),\ldots,\Gr_1(V),*)
\]
Thus, we may view $\phi$ as an non type-preserving automorphism of
the building $\Delta(V)$. Such an isomorphism is called a \emph{duality}.
If $\phi^2=\id_{\Delta(V)}$
(this makes sense in view of the identification
above), then $\phi$ is called a \emph{polarity}.
In this section we study polarities and the related geometries, 
\emph{polar spaces}.

\subsection{Forms and polarities}
\label{SecForms}
\begin{quote}
\emph{We study some basic properties of forms (sesquilinear maps)
and their relation with dualities and polarities.}
\end{quote}
In this section, $V$ is a (possibly infinite dimensional) right
vector space over $D$.
We fix an anti-automorphism $\sigma$ of $D$. Using $\sigma$, we make
the dual space $V^\vee$ into a \emph{right} vector space over $D$,
denoted $V^\sigma$, by defining
$\lambda a=[\xymatrix@1{v\ar@{|->}[r]&a^\sigma\lambda(v)}]$,
for $v\in V$, $\lambda\in V^\vee$
and $a\in D$. We put
\[
\Form_\sigma(V)=\Hom_D(V,V^\sigma).
\]
Given an element $f\in\Form_\sigma(V)$, we write
$f(u,v)=f(u)(v)$; the map $\xymatrix@1{(u,v)\ar@{|->}[r]&f(u,v)}$
is \emph{biadditive}
($\ZZ$-linear in each argument) and \emph{$\sigma$-sesquilinear},
$f(ua,vb)=a^\sigma f(u,v) b$. Note that $\sigma$ is uniquely determined
by the map $f$, provided that $f\neq 0$.
Suppose that that
\[
\xymatrix{{F:V}\ar[r]&{V'}}
\]
is a linear map of vector spaces over $D$. We define
\[
\xymatrix{{\Form_\sigma(V)}\ar@{<-}[r]^{F^*}&{\Form_\sigma(V')}}
\]
by $F^*(f')(u,v)=f'(F(u),F(v))$. 
A similar construction works if $\xymatrix@1{V\ar[r]^F&{V'}}$
is $\theta$-semilinear
relative to a skew field isomorphism
$\xymatrix@1{D\ar[r]^\theta_\cong&{D'}}$; here, we define
$F^*(f')(u,v)=f'(F(u),F(v))^{\theta^{-1}}$ to obtain
\[
\xymatrix{{\Form_\sigma(V)}\ar@{<-}[r]^{F^*}&{\Form_{\sigma'}(V'),}}
\]
where $\sigma=\theta\sigma'\theta^{-1}$.
The group $\GL(V)$ acts thus in a natural way from the right on
$\Form_\sigma(V)$, by putting 
\[
fg=g^*f.
\]
Forms in the same $\GL(V)$-orbit are called \emph{equivalent}; the stabilizer
of a form $f$ is denoted 
\[
\GL(V)_f=\{g\in\GL(V)|\ f(g(u),g(v))=f(u,v)\text{ for all }u,v\in V\}.
\]
The assignment $\xymatrix@1{V\ar@{|->}[r]&{V^\sigma}}$ is a natural
cofunctor on $\Mod_D^\fin$ (or $\Mod_D$, if we allow infinite dimensional
vector spaces) which we also denote by $\sigma$,
\[
\xymatrix{
V \ar[rr]^F &\ar@{->}[dd]^{\sigma}& {V'} \\
\\
{V^\sigma}  && \ar[ll]^{F^\sigma}{{V'}^\sigma}.
}
\]
Let $\xymatrix{{V^\vee}\ar@{<-}[r]^{F^\vee}&{{V'}^\vee}}$ be the dual
or \emph{adjoint} of $\xymatrix{V\ar[r]^F&{V'}}$
(i.e. $F^\vee(\lambda)=\lambda F$). Then set-theoretically, $F^\sigma=F^\vee$.
There is a canonical 
linear injection 
\[
\xymatrix{V\ar[r]^{\mathsf{can}}&{V^{\sigma\sigma}}}
\]
which sends
$v\in V$ to the map
$\mathsf{can}(v)=[\xymatrix@1{\lambda\ar@{|->}[r]&\lambda(v)^{\sigma^{-1}}}]$.
If the dimension of $V$ is finite, then $\mathsf{can}$ is
an isomorphism
(and the data $\sigma$ and $\mathsf{can}$ make the abelian
category $\Mod_D^\fin$
into a \emph{hermitian category} $\mathfrak{Herm}_{D,\sigma}^\fin$).
Given a form $f$, we have a diagram
\[
\xymatrix{
V \ar[rr]^f\ar[ddrr]^{\mathsf{can}} \ar[dd]_{f^\sigma\mathsf{can}}
&& {V^\sigma} \\
\\
{V^\sigma}  && \ar[ll]_{f^\sigma}{{V^\sigma}^\sigma}.
}
\]
Note that
$f^\sigma\mathsf{can}(v)=[\xymatrix@1{u\ar@{|->}[r]&{f(u,v)^{\sigma^{-1}}}}]$.
If both $f$ and $f^\sigma\mathsf{can}$ are injective, then we call the
form $f$ \emph{non-degenerate} (if $\dim(V)$ is finite, then it suffices to
require that $f$ is injective).
For a subspace $U\subseteq V$, put
\begin{align*}
U^{\perp_f}&=\bigcap\{\kernel (f(u))|\ u\in U\}
=\{v\in V|\ f(u,v)=0\text{ for all }u\in U\} \\
{}^{\perp_f}U&=\bigcap\{\kernel (f^\sigma\mathsf{can}(u))|\ u\in U\}
=\{v\in V|\ f(v,u)=0\text{ for all }u\in U\}.
\end{align*}
Thus $f$ is non-degenerate if and only if $V^{\perp_f}=0={}^{\perp_f}V$.

In the finite dimensional case, a non-degenerate form defines a
duality of $\PG(V)$ (by $\xymatrix@1{U\ar@{|->}[r]&{U^{\perp_f}}}$, and
by the Fundamental Theorem of Projective Geometry \ref{FTPG1},
every duality of $\PG(V)$ is obtained in this way).
Also, $\perp_f$ determines the anti-automorphism $\sigma$ 
up to conjugation with
elements of $\Int(D)$. The form $f$ itself is, however, not 
determined by $\perp_f$.
Therefore, we introduce another equivalence relation on forms.
If $f$ is $\sigma$-sesquilinear, and if $s\in D^\times$, then
$\xymatrix@1{sf:(u,v)\ar@{|->}[r]&sf(u,v)}$ is
$\sigma s^{-1}$-sesquilinear,
\[
sf(ua,vb)=sa^\sigma f(u,v)b=(sa^\sigma s^{-1})sf(u,v)b=
a^{\sigma s^{-1}} (sf)(u,v)b.
\]
The forms $f$ and $sf$ are called \emph{proportional}, and we say
that $sf$ is obtained from $f$ by \emph{scaling} with $s$; proportional
forms induce the same dualities. Scaling with $s$ yields an isomorphism
$\xymatrix@1{{\Form_\sigma(V)}\ar[r]^\cong&{\Form_{\sigma s^{-1}}(V)}}$.
\begin{Prop}
Suppose that $\dim(V)$ is finite.
Let $f\in\Form_\sigma(V)$ and $f'\in\Form_{\sigma'}(V)$ be
non-degenerate forms. If $f$ and $f'$ induce the same duality, then
$\sigma'\sigma^{-1}\in\Int(D)$, and $f$ and $f'$ are proportional.
There is a 1-1 correspondence
\[
\xymatrix{
*++{\left\{\txt{Dualities in\\ $\PG(V)$}\right\}}
\ar@{<~>}[rr]&&
*++{\left\{\txt{Proportionality classes of\\
non-degenerate forms}\right\}.}
}
\]
\qed
\end{Prop}
A non-degenerate sesquilinear form which induces a duality has the
property that $(U^{\perp_f})^{\perp_f}=U$ holds for all subspaces $U$.
This condition makes also sense in the infinite dimensional case
if we restrict it to finite dimensional subspaces
(although there, no dualities exist), and
boils down to $f$ being \emph{reflexive};
\[
f(u,v)=0\quad\Longleftrightarrow\quad f(v,u)=0.
\]
In other words, $U^{\perp_f}={}^{\perp_f}U$ holds
for all subspaces $U\subseteq V$.
If $f\neq 0$ is reflexive, then there exists a unique
element $\eps\in D^\times$
such that $f(u,v)=f(v,u)^\sigma\eps$ for all $u,v\in V$,
see Dieudonn\'e \cite{Dieu} Ch.~I \S6.
Furthermore, this implies that $\eps^\sigma=\eps^{-1}$ (choose $u,v$ with
$f(u,v)=1$, then $f(v,u)=\eps$), and $a^{\sigma^2}=a^{\eps^{-1}}$
for all $a\in D$ (consider $f(u,va)$).
A form $h$ which satisfies the identity
\[
h(u,v)=h(v,u)^\sigma\eps\quad\text{ for all }u,v\in V
\]
is called \emph{$(\sigma,\eps)$-hermitian}. The collection of all
$(\sigma,\eps)$-hermitian forms  is a subgroup of $\Form_\sigma(V)$
which we denote $\Herm_{\sigma,\eps}(V)$.
\emph{A non-degenerate $(\sigma,\eps)$-hermitian
form $h$ induces thus an involution $\perp_h$
on the building $\Delta(V)$, for finite dimensional $V$.}
In the finite dimensional setting, we have thus
a 1-1 correspondence
\[
\xymatrix{
*++{\left\{\txt{Polarities in\\ $\PG(V)$}\right\}}
\ar@{<~>}[rr]&&
*++{\left\{\txt{Proportionality classes\\of
non-degenerate\\$(\sigma,\eps)$-hermitian forms}\right\}}
}
\]
If a $(\sigma,\eps)$-hermitian form $h$
is scaled with $s\in D^\times$, then the resulting form $sh$ is
$(\sigma s^{-1},ss^\sigma\eps)$-hermitian.

\subsection{Polar spaces}
\begin{quote}
\emph{We introduce polar spaces as certain point-line geometries.}
\end{quote}
Assume that $h$ is non-degenerate $(\sigma,\eps)$-hermitian.
An element $U\in\Gr_k(V)$ is called
\emph{absolute} if it is incident with its image
$U^{\perp_h}$. If $2k\leq \dim(V)$, then this means that
$U\subseteq U^{\perp_h}$, or, in other words, that $h|_{U\times U}=0$.
A subspace with this property is called \emph{totally isotropic}
(with respect to $h$). Let $\Gr_k^h(V)$ denote the collection of
all totally isotropic $k$-dimensional
subspaces. The maximum number $k$ for which
this set is non-empty is called the \emph{Witt index} $\ind(h)$ of $h$
(If the vector space dimension is infinite, then $\ind(h)$ can of course
be infinite). 
This makes sense also for possibly degenerate $(\sigma,\eps)$-hermitian
forms:
\[
\ind(h)=\max\{\dim(U)|\ U \subseteq V,\ h|_{U\times U}=0\}
\]
If $h$ is non-degenerate, then $2\,\ind(V)\leq \dim(V)$.
Suppose that $h$ is non-degenerate and of finite index $m\geq 2$.
Let $\PG^h(V)$ denote the structure
\[
\PG^h(V)=(\Gr_1^h(V),\ldots,\Gr_m^h(V),*),
\]
and let $\PG^h(V)_{1,2}=(\Gr_1^h(V),\Gr_2^h(V),*)$ denote the corresponding
point-line geometry. This is an example of a \emph{polar space}.
This geometry has one crucial property: given a point $p\in\Gr_1^h(V)$ and
a line $L\in\Gr^h_2(V)$ which are not incident, $p\not\subseteq L$,
there are two possibilities: either there exists an element
$H\in\Gr_3^h(V)$ containing both $p$ and $L$, or there exists precisely
one point $q$ incident with $L$, such that $p\oplus q\in\Gr^h_2(V)$.
Algebraically, this means that either $L\subseteq p^{\perp_h}$ or
$q=L\cap p^{\perp_h}$ (note that $p^{\perp_h}$ is a hyperplane,
so the intersection has at least dimension 1).

A point-line geometry $(\cP,\cL,*)$ is called a (non-degenerate)
\emph{polar space} if it satisfies the following properties.
\begin{description}
\item[PS1]
There exist two distinct lines.
Every line is incident with at least $3$ points.
Two lines which have more than one point in common are equal.
\item[PS2]
Given $p\in\cP$ there exists $q\in\cP$ such that $p$ and $q$ are not joined by
a line.
\item[PS3]
Given a point $p\in\cP$ and a line $\ell\in\cL$, either
$p$ is collinear with every point which is incident with $\ell$, or
with precisely one point which is incident with $\ell$.
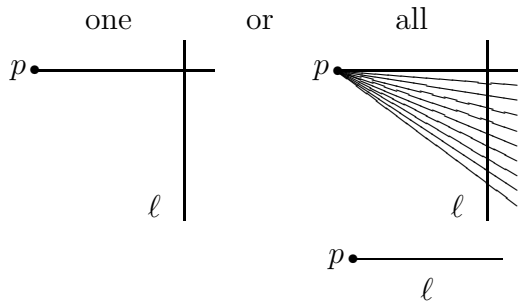
\begin{center}
\begin{tabular}{ccc}
one & or & all \\
\begin{xy} <2cm,0cm>:
(0,1)*@{*}="p"*+!R{p},
(.8,.1)*+U{\ell},
(1,1)*{}="q",
(1.2,1)*{}="a",
(1,1.2)*{}="b",
(1,0)*{}="c",
"p";"a"**\dir{-},
"b";"c"**\dir{-}
\end{xy} &&
\begin{xy} <2cm,0cm>:
(0,1)*@{*}="p"*+!R{p},
(.8,.1)*+U{\ell},
(1,1)*{}="q",
(1.2,1)*{}="a",
(1,1.2)*{}="b",
(1.2,0.9)*{}="q1",
(1.2,0.8)*{}="q2",
(1.2,0.7)*{}="q3",
(1.2,0.6)*{}="q4",
(1.2,0.4)*{}="q5",
(1.2,0.3)*{}="q6",
(1.2,0.2)*{}="q7",
(1.2,0.1)*{}="q8",
(1.2,0.5)*{}="q9",
(1,0)*{}="c",
"p";"a"**\dir{-},
"b";"c"**\dir{-},
"p";"q1"**\dir{-},
"p";"q2"**\dir{-},
"p";"q3"**\dir{-},
"p";"q4"**\dir{-},
"p";"q5"**\dir{-},
"p";"q6"**\dir{-},
"p";"q7"**\dir{-},
"p";"q9"**\dir{-},
"p";"q8"**\dir{-}
\end{xy} \\
&&
\begin{xy} <2cm,0cm>:
(0,0)*@{*}="p"*+!R{p},
(.5,-.2)*+U{\ell},
(1,0)*{}="q",
"p";"q"**\dir{-}
\end{xy}
\end{tabular}
\end{center}
\end{description}
A \emph{subspace} of a geometry satisfying PS1--PS3 is a set $X$
of points with the following
two properties: given two distinct points $p,q\in X$, there exists a
line $\ell$ incident with $p,q$, and if $r$ is also incident with
$\ell$, then $r\in X$.
It is a (non-trivial) fact that every subspace
which contains $3$ non-collinear points is a projective space.
We define the \emph{rank} of $(\cP,\cL,*)$ as
the maximum of the ranks of the subspaces minus 1. 
\begin{description}
\item[PS4]
The rank $m$ is finite.
\item[PS5]
Every subspace of rank $m-2$ is contained in at least $3$ subspaces of
rank $m-1$.
\end{description}
A structure satisfying the axioms PS1--PS5 is called a 
\emph{thick polar space}.
It is easy to see that the structure
\[
\PG^h(V)=(\Gr_1^h(V),\ldots,\Gr_m^h(V),*),
\]
satisfies PS1--PS4, provided that $h$ is non-degenerate and $\ind(h)=m\geq 2$.
Axiom PS5 is more subtle, but is is easy to see that every subspace of
rank $m-1$ is contained in at least $2$ subspaces of rank $m$.
We call such a structure a \emph{weak polar space}
(with thick lines). Our set of axioms is a variation of the
Buekenhout-Shult axiomatization of polar spaces given in
Buekenhout \& Shult \cite{BueShu}.
Similarly as Veblen's axiom in the definition of a
projective geometry in Section \ref{ProGeo},
the one axiom which is geometrically important is the Buekenhout-Shult
'one or all' axiom PS3. By the Fundamental Theorem of Polar Spaces
\ref{FTPS2}, this simple axiom encodes the whole body of geometric
algebra!

We mention some examples of (weak) polar spaces which do not involve hermitian
forms.

\begin{Num}\textbf{Examples}

(0) Let $X,Y$ be disjoint sets of cardinality at least $3$, put
$\cP=X\times Y$, and $\cL=X\cup Y$. By definition, a point $(x,y)$ is incident
with the lines $x$ and $y$. The resulting geometry is a 
weak polar space of
rank $2$; every point is incident with precisely $2$ lines.
This in fact an example of a weak generalized
quadrangle, see Van Maldeghem's article \cite{HVMArt} in these proceedings.

(1) Every (thick) generalized quadrangle (see Van Maldeghem's article
\cite{HVMArt}) is a thick polar space of rank $2$.

(2) Let $V$ be a $4$-dimensional vector space over a skew field $D$, put
$\cP=\Gr_2(V)$ and $\cL=\{(p,A)\in\Gr_1(V)\times\Gr_3(V)|\
p\subseteq A\}$. The incidence is the natural
one (inclusion of subspaces). The resulting polar space which we denote
by $\mathsf{A}_{3,2}(D)$ has rank $3$;
the planes of the polar space are the points and planes of $\PG(V)$,
and every line of this polar space is incident with precisely two planes,
so $\mathsf{A}_{3,2}(D)$ is a weak polar space.
\end{Num}
The next theorem is an important step in the classification of polar
spaces of higher rank (the full classification will be stated in
Section \ref{PSClass}).
\begin{Thm}[Tits]\label{ResiduesAreMoufang}
Let $X$ be a subspace of a polar space $(\cP,\cL,*)$. If $X$ contains
$3$ non-collinear points, then $X$, together with the set of all lines
which intersect $X$ in more than one point, is a self-dual
projective space.
This projective space is either a Moufang plane over some alternative
field $A$, or a desarguesian projective space over some skew field $D$.

\medskip\noindent
\emph{For a proof see Tits \cite{TitsLNM} 7.9, 7.10, and 7.11}.
\qed
\end{Thm}

\subsection{Hermitian forms}
\begin{quote}
\emph{We continue to study properties of hermitian forms.}
\end{quote}
We fix the following data: $D$ is a skew field, $\sigma$ is an
anti-automorphism of $D$, and $\eps\in D^\times$ is an element with
$\eps^\sigma\eps=1$ and $x^{\sigma^2}=x^{\eps^{-1}}$ as in Section
\ref{SecForms}, and $\Herm_{\sigma,\eps}(V)$ is the group of all
$(\sigma,\eps)$-hermitian forms on $V$.
We define the set of \emph{$(\sigma,\eps)$-traces} as
\[
D_{\sigma,\eps}=\{c+c^\sigma\eps|\ c\in D\}
\]
This is an additive subgroup of $(D,+)$, with the property that
$c^\sigma D_{\sigma,\eps}c\subseteq D_{\sigma,\eps}$, for all $c\in D$.
A form $h\in\Herm_{\sigma,\eps}(V)$ is called \emph{trace valued} if
$h(v,v)\in D_{\sigma,\eps}$ holds for all $v\in V$.
We call such a form $h$ \emph{trace $(\sigma,\eps)$-hermitian}.
(If $\mathsf{char}(D)\neq 2$, then it's easy to show that
every hermitian form is trace hermitian.) A hermitian form $h$ is
trace hermitian if and only if it can be written as
\[
h(u,v)=f(u,v)+f(v,u)^\sigma\eps,
\]
for some $f\in\Form_\sigma(V)$.
We denote the group of all trace hermitian
forms by $\mathsf{TrHerm}_{\sigma,\eps}(V)$.
\begin{Lem}
Let $V_0=\{v\in V|\ h(v,v)\in D_{\sigma,\eps}\}$. Then $V_0$ is a
subspace of $V$ containing all totally isotropic subspaces of $V$.
\qed
\end{Lem}
The form induced by $h$ on $V_0/(V_0\cap V_0^{\perp_h})$ is thus
trace hermitian and non-degenerate.
Since we are only interested in the polar space arising from $h$,
we can thus safely assume that $h$ is trace hermitian.

Suppose now that we scale the form with an element $s\in D^\times$.
Let $h'=sh$ and put $\sigma'=\sigma s^{-1}$, and $\eps'=ss^{-\sigma}\eps$.
Then $D_{\sigma',\eps'}=sD_{\sigma,\eps}$. Thus we can achieve that either
$D_{\sigma,\eps}=0$, or that $1\in D_{\sigma,\eps}$. In the first case
we have $\eps'=-1$ and $\sigma=\sigma'=\id_D$ (and then $D$ is commutative),
and in the second case $\eps=1$ and ${\sigma'}^2=\id_D$.
The study of (non-degenerate) trace hermitian forms is thus --- by means of
scaling --- reduced to the following subcases:
\begin{description}
\item[Trace $\sigma$-hermitian forms]\ 

$\sigma^2=\id_D$ (here $\sigma=\id_D$ is allowed if $D$ is commutative)
and $h(u,v)=h(v,u)^\sigma$ for all $u,v\in V$.

If $\sigma=\id_D$, then we call $h$ \emph{symmetric}, and 
(in the non-degenerate case)
$\GL(V)_h=\mathsf{O}(V,h)$ is the  \emph{orthogonal group} of $h$; 
if $\sigma\neq\id_D$ (and if $h$ is non-degenerate), then
$\U(V,h)=\GL(V)_h$ is called the \emph{unitary group} of $h$.
\item[Symplectic forms]\ 

$\sigma=\id_D$ and $h(v,v)=0$  for all $v\in V$ (so $h$ is symplectic).
The group $\GL(V)_h=\Sp(V,h)$ is called the \emph{symplectic group}
(in the non-degenerate case).
\end{description}
Before we consider these groups and forms in more detail, we extend
the whole theory to include (pseudo-)quadratic forms.

\subsection{Pseudo-quadratic forms}
\begin{quote}
\emph{We introduce pseudo-quadratic forms, which are 
certain cosets of ses\-qui\-li\-ne\-ar
forms. In characteristic~$2$, this generalizes trace hermitian forms.}
\end{quote}
In the course of the classification of polar spaces, it
turns out that in characteristic $2$, $(\sigma,\eps)$-trace
hermitian forms are not sufficient; one needs the notion of
a pseudo-quadratic form which is due to Tits. We follow the 
treatment which is now standard and which is based on Bak's concept
\cite{Bak} of form parameters. This is a modification of Tits' original
approach; however, the characteristic $2$ theory of
unitary groups over skew fields is entirely due to Tits, a fact
which is not always properly reflected in books on classical groups
(see e.g. the footnote on p.~190 in Hahn \& O'Meara \cite{HOM}).
Let $\Lambda$ be an additive subgroup of $D$, with
\[
D_{\sigma,-\eps}=\{c-c^\sigma\eps|\ c\in D\}
\subseteq \Lambda\subseteq
D^{\sigma,-\eps}=\{c\in D|\ c^\sigma\eps=-c\},
\]
and with the property that 
\[
s^\sigma\Lambda s\subseteq\Lambda
\]
for all $s\in D$. Such a subset $\Lambda$ is called a \emph{form parameter}.
Given an element $f\in\Form_\sigma(V)$, we define the 
\emph{pseudo-quadratic form} $[f]=(q_f,h_h)$ to be the
pair of maps
\[
\xymatrix @!=1pc {
**[l] {q_f:V} \ar[r] & **[r]{D/\Lambda} &&&&
**[l] {h_f:V\times V}\ar[r] & **[r]D \\
**[l] v \ar@{|->}[r] & **[r] {f(v,v)+\Lambda} &&&&
**[l] {(u,v)} \ar@{|->}[r] & **[r]{f(u,v)+f(v,u)^\sigma\eps .}
}
\]
It is not difficult to see that the map
$\xymatrix@1{f\ar@{|->}[r]&{[f]}}$ is additive;
the kernel is
\[
\Lambda\text{-}\Herm_{\sigma,-\eps}(V)
=\{f\in\mathsf{TrHerm}_{\sigma,-\eps}(V)|\
f(v,v)\in\Lambda\text{ for all }v\in V\}.
\]
The resulting group of \emph{pseudo-quadratic forms}
is denoted $\Lambda$-$\Quad_{\sigma,\eps}(V)$, and we have an exact sequence
\[
\xymatrix{
0 \ar[r]  & {\Lambda\text{-}\Herm_{\sigma,-\eps}(V)} \ar[r] &
{\Form_\sigma(V)} \ar[r]
 & {\Lambda\text{-}\Quad_{\sigma,\eps}(V)} \ar[r] & 0 }
\]
which is compatible with the $\GL(V)$-action on $\Form_\sigma(V)$;
furthermore, scaling is a well-defined process on pseudo-quadratic forms.
Let $\psi(c)=c+c^\sigma\eps$. Then
$\mathsf{ker}(\psi)=D^{\sigma,-\eps}\supseteq\Lambda$;
consequently, there is a well-defined map
$\xymatrix@1{{D/\Lambda}\ar[r]^(.6){\bar\psi}&D}$
sending $c+\Lambda$ to $\psi(c)$,
\[
\xymatrix{
 & *++{\Lambda} \ar@{^{(}->}[d]\ar[dr] \\
0 \ar[r]  & {D^{\sigma,-\eps}} \ar[r] & {D} \ar[r]^{\psi} \ar[dr]
 & {D_{\sigma,\eps}} \ar[r] & 0 .\\
 &&& {D/\Lambda} \ar[u]_{\bar\psi} }
\]
Thus we have $\bar\psi(q_f(v))=h(v,v)$;
in particular,
\[
q_f(v)=0\quad\Longrightarrow\quad h_f(v,v)=0.
\]
We call a pseudo-quadratic form \emph{non-degenerate}
if $h_f$ is non-degenerate (this differs from Tits' notion of
non-degeneracy \cite{TitsLNM}, \cite{BruhTits};
we'll come back to that point later).
Similarly as before, a subspace $U$ is called
\emph{totally isotropic} if $q_f$ and $h_f$ vanish on $U$;
the collection of all $k$-dimensional totally isotropic subspaces
is denoted $\Gr_k^{[f]}(V)$, and the Witt index $\ind[f]$ is defined
in the obvious way. Note that 
\[
\Gr_k^{[f]}(V)\subseteq\Gr_k^{h_f}(V).
\]
Suppose that $V,V'$ are vector spaces over $D$, and that
\[
\xymatrix{{F:V} \ar[r] & {V'}}
\]
is linear. We define a map 
\[
\xymatrix{{\Lambda\text{-}\Quad_{\sigma,\eps}(V)}
\ar@{<-}[r]^{F^*}&{\Lambda\text{-}\Quad_{\sigma,\eps}(V'),}}
\]
by $\xymatrix@1{{[f']}\ar@{|->}[r]&{[F^*(f')]}}$
which we denote also by $F^*$. 
A similar construction works if $\xymatrix@1{V\ar[r]^F&{V'}}$
is $\theta$-semilinear
relative to $\xymatrix@1{D\ar[r]^\theta_\cong&{D'}}$.
The group
\[
\U([f])=\GL(V)_{[f]}
\]
is the group of all isometries of $(V,[f])$. Let $g\in\GammaL(V)$.
If there exists an element $s\in D$ such that $g^*[f]=[sf]$, then
$g$ is called a \emph{semi-similitude}; if $g$ is linear, then it is called a
\emph{similitude}. The corresponding groups are denoted
\[
\xymatrix{{\U([f])}\ar@{^{(}->}[r]^{\unlhd}&
{\GU([f])}\ar@{^{(}->}[r]^{\unlhd}&{\GammaU([f]).}}
\]
Not every automorphism $\theta$ of $D$ can appear in $\GammaU([f])$;
a necessary and sufficient condition is that
\[
[\theta,\sigma]\in\Int(D);
\]
up to inner automorphisms, $\theta$ has to centralize $\sigma$.
\begin{Num}\textbf{The case when $h_f$ is degenerate}

\noindent
There is one issue which we have to address. In characteristic $2$, it
is possible that $V^{\perp_{h_f}}\neq 0$, while
$q_f^{-1}(0)\cap V^{\perp_{h_f}}=0$. Let's call such a pseudo-quadratic
form  \emph{slightly degenerate}. This case can be reduced to the
non-degenerate case as follows.
Let $V'=V/V^{\perp_{h_f}}$, and let 
\[
\Lambda'=\{c\in D|\ c+\Lambda\in q_f(V^{\perp_{h_f}})\}.
\]
It can be checked that $\Lambda'$ is a form parameter.
Define $(\tilde q,\tilde h)$ on $V'$ by
\[
\tilde h(u+V^{\perp_{h_f}},v+V^{\perp_{h_f}})=h(u,v)\text{ and }
\tilde q(v+V^{\perp_{h_f}})=f(v,v)+\Lambda'.
\]
One can check that this pair is a pseudo-quadratic form
$[\tilde f]=(\tilde q,\tilde h)$;
there is a canonical bijection
\[
\xymatrix{{\Gr_k^{[f]}(V)}\ar[r]&{\Gr_k^{[\tilde f]}(V').}}
\]
Furthermore, there is a corresponding isomorphism
$\GL(V)_{[f]}\cong\GL(V')_{[\tilde f]}$.
Here is an example. Let $D$ be a perfect field of characteristic $2$,
let $V=D^5$, and let $f$ denote the bilinear form given by the matrix
\[
f\sim
\begin{pmatrix}
0 & 1 \\
 & 0 \\
 && 0 & 1 \\
 &&& 0 \\
 &&&& 1
\end{pmatrix}
\]
Thus $q_f(x)=x_1x_2+x_3x_4+x_5^2$. The associated bilinear form
$h_f$ is symplectic and degenerate; its matrix is
\[
h_f\sim
\begin{pmatrix}
0 & 1 \\
1 & 0 \\
 && 0 & 1 \\
 && 1 & 0 \\
 &&&& 0
\end{pmatrix}
\]
The process above gives us an isomorphism
between an orthogonal generalized quadrangle and a symplectic generalized
quadrangle, and a group isomorphism
\begin{align*}
\mathsf{O}(q_f,D) & \cong\mathsf{Sp}_4(D)\\
\mathsf{O}(5,2^k) & \cong\mathsf{Sp}(4,2^k) \text{ for the finite
case }D=\FF_{2^k}.
\end{align*}
\end{Num}
In Tits \cite{TitsLNM} and Bruhat-Tits \cite{BruhTits},
the chosen form parameter is always
the minimal one, $\Lambda=D_{\sigma,-\eps}$. Therefore, Tits
allows his forms to be slightly degenerate (in our terminology).
The resulting theory is the same; the choice of a bigger $\Lambda$
makes the vector spaces smaller and avoids degenerate hermitian forms,
which is certainly more elegant;
the expense is that in this way we don't really see groups like
$\mathsf{O}(5,2^k)$, since they are identified with their isomorphic
images belonging to non-degenerate forms,
$\mathsf{O}(5,2^k)\cong\mathsf{Sp}(4,2^k)$ ---
one should keep that in mind.

\subsection{Properties of form parameters}
\begin{quote}
\emph{We discuss some properties of pseudo-quadratic forms and their
form parameters.}
\end{quote}
In general, we have
\[
\Lambda=D^{\sigma,-\eps}\quad\Longrightarrow\quad
\biggl( q_f(v,v)=0\quad\Longleftrightarrow\quad h_f(v,v)=0\biggr)
\]
so the theory of trace hermitian forms is included in the pseudo-quadratic
forms as a subcase. So the question is:
\begin{center}
\emph{Why pseudo-quadratic forms?}
\end{center}
The answer is given by the Fundamental Theorem of Polar Spaces \ref{FTPS2}.
But first, we mention a few cases where pseudo-quadratic forms are \emph{not}
necessary. As we mentioned above, this is the case if
$\Lambda=D^{\sigma,-\eps}$. Now if $\mathsf{char}(D)\neq 2$, then
$D_{\sigma,-\eps}=D^{\sigma,-\eps}$, so in characteristic different
from $2$, $(\sigma,\eps)$-hermitian forms suffice.
\begin{Lem}
If $\mathsf{char}(D)\neq 2$, then there is a natural isomorphism
\[
\mathsf{TrHerm}_{\sigma,\eps}(V)\cong
\Lambda\text{-}\Quad_{\sigma,\eps}(V).
\]
\end{Lem}
We consider some more special cases of form parameters.
Note that
\[
\Lambda=0\quad\Longrightarrow\quad\biggl(\eps=1\text{ and }\sigma=\id_D
\text{ and $D$ commutative}\biggr).
\]
Suppose now that $D$ is commutative. If $\sigma=\id_D$, then
$D_{\id_D,-\eps}=D(1-\eps)$. So either $\Lambda=D$, or $\eps=1$.
If $D\neq \Lambda\neq 0$, then we have necessarily $\mathsf{char}(D)=2$,
and $\Lambda$ is a $D^2$-submodule of $D$, and $D$ is not perfect.
\begin{Lem}
Suppose that $D$ is commutative and that $\sigma=\id_D$.
If $0\neq \Lambda\neq D$, then $\mathsf{char}(D)=2$, the field
$D$ is not perfect, and $\Lambda$ is a $D^2$-submodule of $D$. 
If $\Lambda=D$, then $\eps=-1$.
\qed
\end{Lem}
Suppose now that $D$ is commutative and that $\sigma\neq \id_D$.
Then $\sigma^2=\id_D$, and $D_{\sigma,-\eps}\neq 0$. Let $K\subseteq D$
denote the fixed field of $\sigma$.
After scaling, we may assume that $1\in D_{\sigma,-\eps}$, which implies that
$\eps=-1$. Then $D^{\sigma,1}=K$, and if $\mathsf{char}(D)\neq 2$, then
$\Lambda=K$.
For $\mathsf{char}(D)=2$ we put $\psi(c)=c+c^\sigma$; then we
have an exact sequence
\[
\xymatrix{
0 \ar[r]  & K \ar[r] &
{D} \ar[r]^\psi & {D_{\sigma,1}} \ar[r] & 0 }
\]
of finite dimensional vector spaces over $K$, so $D_{\sigma,1}=K$,
regardless of the characteristic.
\begin{Lem}
\label{FormParametersInFields}
Suppose that $D$ is commutative and that $\sigma\neq\id_D$.
Then $\sigma^2=\id_D$, and the form parameters are the left translates
$sK$ of the fixed field $K$ of $\sigma$, for $s\in D^\times$; in
particular, $\Lambda=D^{\sigma,\eps}$.
\qed
\end{Lem}
\begin{Cor}
Form parameters and pseudo-quadratic forms over perfect fields
(in particular, over finite fields) are not important.
\qed
\end{Cor}
This explains why form parameters and proper pseudo-quadratic forms
($\sigma\neq\id_D$) are not an issue in finite geometry, and why
they don't appear in books on algebraic groups over algebraically
closed fields.
Note also that if $D$ is algebraically closed and if $\sigma\neq 1$,
then the fixed field of $\sigma$ is a real closed field.

The next result is less obvious and was pointed out to me by 
Richard Weiss.
\begin{Prop}[Finite form parameters]
Suppose that $\Lambda\neq 0$ is a finite form parameter. Then $D$ is finite.

\proof
If $\sigma=\id_D$, then $D$ is commutative. If $\mathsf{char}(D)=2$, then
$\Lambda\subseteq D$ is a $D^2$-module.
If $\mathsf{char}(D)\neq 2$, then $\Lambda=D$;
in any case, $\Lambda$ has a subset of the same cardinality as $D$.

If $\sigma\neq\id_D$, then it suffices to consider the minimal case where
$\Lambda=D_{\sigma,-\eps}$. We rescale in such a way
that $1\in D_{\sigma,-\eps}$; then $\sigma$ is an involution and
$\eps=-1$. Put $\psi(x)=x+x^\sigma$ and let
$\lambda=\psi(x)\in\Lambda$. We claim that
$\lambda^k\in\Lambda$, for all $k\geq 1$. 
Indeed, $(x+x^\sigma)^k$ can be written as a sum of $2^k$ monomials
of the form 
\[
x^{(\sigma^{\nu_1})}x^{(\sigma^{\nu_2})}\cdots x^{(\sigma^{\nu_k})},
\]
where $\nu_i\in\{0,1\}$. Now
\[
x^{(\sigma^{\nu_1})}x^{(\sigma^{\nu_2})}\cdots x^{(\sigma^{\nu_k})}+
x^{(\sigma^{\nu_k+1})}x^{(\sigma^{\nu_{k-1}+1})}\cdots x^{(\sigma^{\nu_1+1})}
=
\psi(x^{(\sigma^{\nu_1})}x^{(\sigma^{\nu_2})}\cdots x^{(\sigma^{\nu_k})}).
\]
If $\Lambda$ is finite, then all elements $\lambda\in\Lambda\setminus\{0\}$
have
thus finite multiplicative order. This implies that $D$
is commutative, see Herstein \cite{Herstein} 
Cor.~2 p.~116. From Lemma \ref{FormParametersInFields} above, we see that
$\mathsf{card}(\Lambda)+\mathsf{card}(\Lambda)=\mathsf{card}(D)$,
so $D$ is finite.
\qed
\end{Prop}
(It is a consequence of this proposition that there exist no
semi-finite spherical irreducible Moufang buildings: if one panel of such
a building is finite, then every panel is finite. This is a problem
about Moufang polygons, and, by the classification due to Tits and Weiss
(as described by
Van Maldeghem in these proceedings \cite{HVMArt}, the only difficult
case is presented by
the classical Moufang quadrangles associated to hyperbolic
spaces of rank $2$; there, the line pencils are parametrized by the set
$\Lambda$.)

This section shows that pseudo-quadratic forms are important only if
either $D$ is a finite field of characteristic $2$ and if $\sigma=\id_D$
(and then we are dealing with quadratic forms),
or if $D$ is a non-commutative skew field of characteristic $2$.
On the other hand, none of the results about classical groups becomes
really simpler if pseudo-quadratic forms are excluded, so we stick with them.

\subsection{Polar spaces and pseudo-quadratic forms}
\label{PSClass}
\begin{quote}
\emph{We state the classification of polar spaces.}
\end{quote}
Suppose that $[f]$ is a non-degenerate pseudo-quadratic form of
(finite) index $\ind[f]=m\geq 2$. Let
$\PG^{[f]}(V)=(\Gr_1^{[f]}(V),\ldots,\Gr_m^{[f]}(V),*)$ and
\[
\PG^{[f]}(V)_{1,2}=(\Gr_1^{[f]}(V),\Gr_2^{[f]}(V),*).
\]
Then $\PG^{[f]}(V)_{1,2}$ is a (possibly weak) polar space and a
subgeometry of
the (possibly weak)
polar space $\PG^{h_f}(V)_{1,2}$. A polar space isomorphic
to such a space is called \emph{embeddable}.
Here is the first analogue of the Fundamental Theorem of Projective
Geometry.
\begin{Thm}[Fundamental Theorem of Polar Spaces, I]\

\noindent
Let
\[
\xymatrix{{\PG^{[f]}(V)_{1,2}}\ar[r]_\cong^\phi&{\PG^{[f']}(V')_{1,2}}}
\]
be an isomorphism of embeddable (weak) polar spaces of finite ranks
$m,m'\geq 3$.
Then $m=m'$, and
there exists an isomorphism of skew fields
$\xymatrix@1{D\ar[r]^\theta_\cong&{D'}}$
and a $\theta$-semilinear isomorphism
$\xymatrix@1{V\ar[r]^\Phi_\cong&{V'}}$ such that
the pull-back $\Phi^*([f'])$ is proportional to $[f]$.

\medskip\noindent\em
For a proof see Tits \cite{TitsLNM} Ch.~8, or Hahn \& O'Meara 8.1.5.
\end{Thm}
This result and the next one are partly due to Veldkamp;
the full results were
proved by Tits. Cohen \cite{Cohen} and Scharlau \cite{RScharlau}
are good references for the classification, and for newer results
in this area. The situation is more complicated for embeddable polar
spaces of rank $2$; we refer to Tits \cite{TitsLNM} Ch.~8.
The theorem above deals with embeddable polar spaces. In fact,
all polar spaces of higher rank are embeddable.

\begin{Thm}[Fundamental Theorem of Polar spaces, II]\ 

\noindent
Suppose that $(\cP,\cL,*)$ is a (weak) polar space of rank $m\geq 4$. Then
$(\cP,\cL,*)$ is embeddable.

\medskip\noindent\em
For a proof see Tits \cite{TitsLNM} 8.21, combined with
Thm.~\ref{ResiduesAreMoufang}. See also Scharlau \cite{RScharlau}
Sec.~7.
\end{Thm}
This result is not true for polar spaces of rank $3$; there exist
polar spaces which have Moufang planes over alternative fields
as subspaces, and such a polar space cannot be embeddable.
However, this is essentially the only thing which can go wrong.
\begin{Thm}[Fundamental Theorem of Polar spaces of rank $3$]
\label{FTPS2}
\

\noindent
Let $(\cP,\cL,*)$ be a polar space of rank $3$. If $(\cP,\cL,*)$ is not
embeddable, then either there exists a proper alternative field
$A$, and the maximal subspace are projective Moufang planes over $A$,
or $(\cP,\cL,*)\cong\mathsf{A}_{3,2}(D)$, for some skew field $D$.

\medskip\noindent\em
For a proof see Tits \cite{TitsLNM} 7.13, p.~176, and 9.1, combined with
Thm.~\ref{ResiduesAreMoufang}. See also Scharlau \cite{RScharlau}
Sec.~7.
\qed
\end{Thm}
The polar spaces containing proper Moufang planes are related to
exceptional algebraic groups of type $E_7$; these are the only polar
spaces of higher rank which do not come from classical groups.
If $D$ is commutative, then $\mathsf{A}_{3,2}(D)$ is related to the
Klein correspondence,  $\mathbf{D}_3=\mathbf{A}_3$.

Finally, we should mention the following result which is a consequence
of Tits' classification.
\begin{Prop}[Tits]
Let $(\cP,\cL,*)$ be a weak polar space of rank $m\geq 3$, such that
every subspace of rank $m-2$ is incident with precisely two subspaces of rank
$m-1$. Then either $(\cP,\cL,*)\cong\mathsf{A}_{3,2}(D)$, or
$(\cP,\cL,*)\cong\PG^h(V)_{1,2}$, where $V$ is a $2m$-dimensional vector
space over a field $D$, and $h$ is a non-degenerate symmetric bilinear form
(i.e. $\sigma=\id_D$ and $\eps=1$) of index $m$ (in other words, $V$ is
a hyperbolic module of orthogonal type, see the next section).
\qed
\end{Prop}
Thick polar spaces of rank $2$ are the same as generalized quadrangles.
Similarly as projective planes, these geometries can be rather 'wild'
and there is no way to classify them. The classification of the
Moufang quadrangles due to Tits and Weiss is  a major
milestone in incidence geometry. For results about generalized quadrangles
we refer to Van Maldeghem's article \cite{HVMArt} in these proceedings.
Polar spaces of possibly infinite rank where considered by Johnson \cite{PJ}.

\subsection{Polar frames and hyperbolic modules}
\begin{quote}
\emph{We show how the classification of pseudo-quadratic forms is
reduced to the anisotropic case.}
\end{quote}
Let $V$ be an $m$-dimensional vector space over $D$,
and put 
\[
H=V\oplus V^\sigma.
\]
We define a form $f\in\Form_\sigma(H)$ by
\[
f((u,\xi),(v,\eta))=\xi(v).
\]
As a matrix, $f$ is represented as
\[
f\sim
\begin{pmatrix}
0 & 1_m  \\
0& 0
\end{pmatrix}
\]
where $1_m$ denotes the $m\times m$ unit matrix.
The space $H$ with the pseudo-quadratic form $[f]$ (relative
to a form parameter $(\Lambda,\sigma,\eps)$) is called
a \emph{hyperbolic module} of rank $m$; if $m=1$ then $H$ is
$2$-dimensional and we call it a \emph{hyperbolic line}.
(Hyperbolic lines are often called \emph{hyperbolic planes}; this
depends on the viewpoint, linear algebra vs. projective geometry.)
The following theorem is crucial.
\begin{Thm}
Let $[f]$ be a non-degenerate pseudo-quadratic form of finite
Witt index $\ind([f])=m$ in a vector space $V$.
Then there exists a hyperbolic module $H$ of rank $m$ in $V$, and
$V$ splits as an orthogonal sum
\[
V=H\oplus V_0,
\]
where $V_0=H^{\perp_{h_f}}$. If $H'\subseteq V$ is another hyperbolic
module of rank $m$, then there exists an isometry of $V$ which maps
$H'$ onto $H$.

\medskip\noindent\em
This follows from Witt's Theorem, see Hahn \& O'Meara 6.1.12,
6.2.12 and 6.2.13 --- the proofs apply despite the fact that
Hahn \& O'Meara work always with finite dimensional vector spaces.
What is needed in their proof is only that $H$ has finite
dimension.
\qed
\end{Thm}
Thus, the subspace $V_0$ is unique up to isometry. This subspace
(together with the restriction of $[f]$) is called the
\emph{anisotropic kernel} of $[f]$. Since the hyperbolic module
has a relatively simple structure, the study of pseudo-quadratic forms
is reduced to the anisotropic case; a pseudo-quadratic form is
determined its Witt index and its anisotropic kernel
(and by $\sigma,\eps,\Lambda$, of course).

The \emph{building associated to a polar space}
$(\cP,\cL,*)$ is constructed as follows.
If $(\cP,\cL,*)$ is thick, then the vertices are the subspaces of the polar
space, and the simplices are sets of pairwise incident vertices. The resulting
building has rank $m$ and type $\mathbf C_m$, see Tits \cite{TitsLNM} Ch.~7. 
If the polar space is weak, then a new geometry is introduced: the vertices
are all subspaces of rank
different from $m-2$, and two vertices are called incident if one contains
the other, or if their intersection has rank $m-2$. This is again an
$m$-sorted structure (there are two classes of subspaces of rank $m-1$),
and the resulting simplicial complex is a building of type $\mathbf D_m$,
see Tits \cite{TitsLNM} Ch.~6 and 7.12, 8.10.
It is easy (but maybe instructive) to check that this makes the weak polar
space $\mathsf{A}_{3,2}(D)$ into the building $\Delta(D^4)$ obtained from
the projective space $\PG(D^4)$.
The buildings related to polar spaces are also discussed in
Brown \cite{Brown}, Cohen \cite{Cohen},
Garrett \cite{Garrett},
Ronan \cite{Ronan}, Scharlau \cite{RScharlau}, and
Taylor \cite{Taylor}.

Finally, we mention the classical groups obtained from non-degenerate
pseudo-quadratic forms of index $m\geq 2$.
Scaling the form by a suitable constant
$s\in D^\times$,
the following cases appear.

\begin{description}
\item[Symplectic groups]\

This is the situation when $(\sigma,\eps,\Lambda)=(\id_D,-1,D)$.
Here $D$ is commutative, $q_f=0$ and $h_f$ is alternating.
The dimension of $V$ is even (and $V$ is hyperbolic),
and $2\,\ind(h_f)=\dim(V)$.
The corresponding polar space is thick.

\item[Orthogonal groups]\

This is the situation when $(\sigma,\eps,\Lambda)=(\id_D,1,0)$.
Here $D$ is commutative and $h_f$ is symmetric,
and $2\,\ind[f]\leq\dim(V)$.
The corresponding polar space is thick if and only if
$2\,\ind[f]<\dim(V)$. If $2\,\ind[f]=\dim(V)$, then the
corresponding building is the $\mathbf D_m$-building (the
\emph{oriflamme geometry}) described above, and $V$ is hyperbolic.

\item[Defective orthogonal groups]\

This is the situation when $(\sigma,\eps)=(\id_D,1)$ and
$0\neq\Lambda\neq D$.
Here $D$ is commutative and $h_f$ is symmetric. This occurs only
in characteristic $2$ over non-perfect fields (in the perfect case,
$\Lambda=D$ and we are in the symplectic case).

\item[Classical unitary groups]\

This is the situation when $\sigma\neq\id_D=\sigma^2$, $\eps=1$
and $\Lambda=D^{\sigma,-1}$.
Here $D$ need not be commutative and $h_f$ is $(\sigma,1)$-hermitian.
Since $\Lambda$ is maximal, the hermitian form $h_f$ describes
$[f]$ completely and $q_f$ is not important.
The corresponding polar space is thick.

\item[Restricted unitary groups]\

This is the situation when $\sigma\neq\id_D=\sigma^2$, $\eps=1$
and $\Lambda<D^{\sigma,-1}$.
Here $D$ is of characteristic $2$ and not commutative.
The corresponding polar space is thick.
\end{description}

\subsection{Omissions}
By now it should be clear that the theory of pseudo-quadratic forms
and the related geometries is rich, interesting, and sometimes difficult.
There are many other interesting topics which we just mention without
further discussion.

\begin{description}
\item[Root elations]\

Root elations in polar spaces are more complicate than root elations
in projective spaces. This is due to the fact that there are two
types of half-apartments and, consequently, two types of root groups.
One kind is isomorphic to the additive group of $D$, while the other
is related to the anisotropic kernel $V_0$ of $V$, and to the form
parameter $\Lambda$. These root groups are nilpotent of class $1$ or $2$.
We refer to Van Maldeghem \cite{HVM} for a detailed description of the
root groups. The root elations are Eichler transformations (also
called Siegel transformations), which are special products of
transvections, and the group generated by these maps is the
elementary unitary group $\mathsf{EU}([f])$.

\item[K-theory]\

Starting with the abstract commutator relations for the root groups of
a given apartment, one can construct unitary version of the Steinberg
groups, and unitary K-groups. Hahn \& O'Meara \cite{HOM} give a
comprehensive introduction to the subject (for hyperbolic $V$).
There is a natural map $\xymatrix@1{{\sK_0(D)}\ar[r]&{\mathsf{KU}_0(D)}}$
whose cokernel
is the Witt group of $D$, another important invariant.
\item[Permutation groups]\

If $\ind[f]\geq 2$, then the action of the unitary group on
$\Gr_1^{[f]}(V)$ is not $2$-transitive. Instead, one obtains
interesting examples of permutation groups of rank $3$
(i.e. with $3$ orbits in $\Gr_1^{[f]}(V)\times\Gr_1^{[f]}(V)$).

\item[Moufang sets]\

If $\ind[f]=1$, then the corresponding unitary group is $2$-transitive
on $\Gr_1^{[f]}(V)$, and there is a natural Moufang set structure.

\item[Isomorphisms]\

As in the linear case, one can ask whether two unitary groups can
be (abstractly or as permutation groups) isomorphic. Indeed, there
are several interesting isomorphisms related to the Klein correspondence
and to Cayley algebras. Many results
in this direction can be found in Hahn \& O'Meara \cite{HOM}.
\end{description}

\end{document}